\documentclass[12pt, reqno]{amsart}

\usepackage[T1]{fontenc}
\usepackage{lmodern}
\usepackage[margin=1in]{geometry}
\usepackage{amsmath,amssymb,amsthm,mathtools}
\usepackage{hyperref}
\usepackage{xcolor}
\usepackage{enumitem}
\usepackage{cleveref}
\usepackage{listings}
\usepackage{fontawesome} 
\usepackage{float}

\makeatletter
\AtBeginDocument{%
	\let\c@lstlisting\relax
	\newcounter{lstlisting}[section]
	\gdef\thelstlisting%
	{\ifnum \c@section>\z@ \thesection.\fi \@arabic\c@lstlisting}}
\makeatother

\definecolor{keywordcolor}{rgb}{0.7, 0.1, 0.1}   
\definecolor{keywordcolor}{cmyk}{0,0.90,0.86,0} 
\definecolor{tacticcolor}{rgb}{0.0, 0.1, 0.6}    
\definecolor{tacticcolor}{cmyk}{1,0.1,0,0.1} 
\definecolor{commentcolor}{rgb}{0.4, 0.4, 0.4}   
\definecolor{commentcolor}{cmyk}{1,0.1,0,0.1} 
\definecolor{symbolcolor}{rgb}{0.0, 0.1, 0.6}    
\definecolor{symbolcolor}{cmyk}{1,0.1,0,0.1} 
\definecolor{sortcolor}{rgb}{0.1, 0.5, 0.1}      
\definecolor{sortcolor}{cmyk}{0.20,0,1,0.19}    
\definecolor{attributecolor}{rgb}{0.7, 0.1, 0.1} 
\definecolor{attributecolor}{cmyk}{0,0.90,0.86,0} 

\let\lean\lstinline

\def\mmlean{\lstinline[mathescape]}

\usepackage{xspace}

\newcommand*{\ud}{\, \mathrm{d}}
\newcommand*{\nat}{\mathbb{N}} 
\newcommand*{\complex}{\mathbb{C}} 
\newcommand*{\real}{\mathbb{R}} 
\newcommand{\nnreal}{\real_{\ge0}}
\newcommand{\ennreal}{\overline{\real_{\ge0}}}

\newcommand{\mathlib}{\texttt{Mathlib}\xspace}

\newcommand{\extlink}{~\ensuremath{{}^\text{\faExternalLink}}}

\theoremstyle{definition}
\newtheorem{definition}{Definition}[section]

\hypersetup{
  bookmarksopen,bookmarksnumbered,
  colorlinks = true,
  linkcolor  = blue,
  citecolor  = red,
  urlcolor   = gray,
  breaklinks
}
\usepackage[backend=biber,maxbibnames=5,maxalphanames=4,style=alphabetic]{biblatex}
\usepackage[autostyle]{csquotes}
\theoremstyle{plain}
\newtheorem{theorem}{Theorem}[subsection]

\DeclareMathOperator{\Lip}{\operatorname{Lip}}

\newcommand{\pc}{{\mathrm{c}}}
\newcommand{\ps}{{\mathrm{s}}}

\newcommand{\fc}{{\Omega}}

\newcommand{\scI}{{\mathcal{I}}} 
\newcommand{\tQ}{{Q}}
\newcommand{\mfa}{{\theta}}
\newcommand{\mfb}{{\zeta}}
\newcommand{\Mf}{{\Theta}}
\newcommand{\fcc}{{\mathcal{Q}}}
\newcommand{\fp}{{\mathfrak p}}
\newcommand{\fP}{{\mathfrak P}}
\newcommand{\fq}{{\mathfrak q}}

\title{\textbf{Formalizing Carleson's Theorem in Lean}}
\author[Becker et al.]{Lars Becker}
\address{Lars Becker,
	Princeton University}
\email{lbecker@math.princeton.edu}

\author[]{Mar\'ia Inés de Frutos-Fernández}
\address{Mar\'ia Inés de Frutos-Fernández,
	University of Bonn}
\email{midff@math.uni-bonn.de}

\author[]{Leo Diedering}
\address{Leo Diedering,
	University of Bonn}
\email{leo.diedering@uni-bonn.de}

\author[]{Floris van Doorn}
\address{Floris van Doorn,
	University of Bonn}
\email{vdoorn@math.uni-bonn.de}

\author[]{Sébastien Gouëzel}
\address{Sébastien Gouëzel,
	University of Rennes}
\email{sebastien.gouezel@univ-rennes1.fr}

\author[]{Evgenia Karunus}
\address{Evgenia Karunus,
  University of Bonn}
\email{lakesare@gmail.com}

\author[]{Edward van de Meent}
\address{Edward van de Meent,
  Utrecht University}
\email{edwardvdmeent@gmail.com}

\author[]{Pietro Monticone}
\address{Pietro Monticone,
Harmonic}
\email{pietro.monticone@harmonic.fun}

\author[]{Jasper Mulder-Sohn}
\address{Jasper Mulder-Sohn,
	The Hague}
\email{jasper.mulder@planet.nl}

\author[]{Jim Portegies}
\address{Jim Portegies,
  Eindhoven University of Technology}
\email{j.w.portegies@tue.nl}

\author[]{Joris Roos}
\address{Joris Roos,
  University of Massachusetts Lowell}
\email{joris\_roos@uml.edu}

\author[]{Michael Rothgang}
\address{Michael Rothgang,
	University of Bonn}
\email{rothgang@math.uni-bonn.de}

\author[]{James Sundstrom}
\address{James Sundstrom,
	Baruch College}
\email{james.sundstrom@baruch.cuny.edu}

\author[]{Jeremy Tan}
\address{Jeremy Tan,
	National University of Singapore}
\email{jtjierui@comp.nus.edu.sg}

\date{\today}

\begin{document}

\begin{abstract}
We present the formalization of Carleson's theorem in the proof assistant \emph{Lean}.
This paper describes the mathematical content, the organization of the project, the blueprint,
and the main design decisions behind the formalization. It is the result of a large
collaborative effort, written and developed in public.
\end{abstract}

\maketitle

\vspace{-\baselineskip}
\tableofcontents

\section{Introduction}\label{sec:intro}

In 1966, Carleson proved that the Fourier series of a square-integrable periodic function
converges to that function at almost every point~\cite{carleson},
settling a conjecture of Luzin~\cite{Luzin13} that had been open for more than fifty years.
In this article, we report on the formalization in the proof assistant \emph{Lean}~\cite{moura2021lean}
of a recent, far-reaching generalization of this theorem to doubling metric measure
spaces~\cite{ThieleCarlesonPreprint}, called the metric space Carleson theorem, and of the deduction of Carleson's statement from it.

This formalization project is notable for a variety of reasons.
The first reason is the proof itself.
Carleson's argument was famously difficult, and later proofs due to
Fefferman~\cite{fefferman} and to Lacey and Thiele~\cite{lacey-thiele},
while conceptually more transparent, remain long and highly technical~\cite[p. 170]{MR3334208}.
The second reason is that this is the first formalization of a state-of-the-art theorem from
harmonic analysis, an area that has not received very much attention from the formalization
community.
The third reason is the highly collaborative nature of this project.
It was a large-scale collaborative formalization,
developed in public and open to anyone who wanted to contribute:
all authors formalized a significant part of the theorems.
The last reason is that this is one of a handful of examples in mathematics where the manual formalization
of a preprint was finished prior to the submission of the preprint to a journal.

\subsection{Brief history and significance of Carleson's theorem}\label{subsec:history}

One of the oldest questions in Fourier analysis consists of understanding the nature of convergence of the partial Fourier sums $S_n f$ of a function $f$ (see Definition~\ref{def:fourier-series}).
There are two main kinds of convergence that are commonly considered. First, one can ask whether the partial Fourier sums $S_n$ converge to the function $f$ in a suitable function
space. There are several well-known answers to this question. For instance, convergence in the Hilbert
space sense for $f$ in $L^2(0,2\pi)$ was established in the first decade of
the twentieth century as a consequence of the rapid development of Lebesgue
integration theory, and convergence in the sense of
distributions for general distributional $f$ was discovered a few decades
after Lebesgue integration. For some other natural spaces, such as
$L^1(0,2\pi)$, there is no guarantee of convergence in the norm of that space,
even if $f$ is in the space.

A second question that turns out to be much more difficult is whether the sequence $\{S_n f(x)\}_{n\in\nat}$ converges pointwise to $f(x)$ for almost every $x$.
Luzin \cite{Luzin13} conjectured that this convergence would hold for functions in the
$L^2(0,2\pi)$ space. In 1923, Kolmogorov found an example \cite{Kolmogorov} of an $L^1$ function whose Fourier series diverges at almost every point, which seemed to suggest that the conjecture might be false. However, when in the 1960s Carleson tried to construct a counterexample to Luzin's conjecture, he realized that a potential counterexample would have to satisfy so many properties that they would in fact be contradictory. He was therefore able to conclude that such a counterexample could not exist, hence proving Luzin's conjecture \cite{carleson}.

Carleson's result was soon after generalized by Hunt \cite{MR238019} to the case of $L^p$ functions with $p>1$ and for functions in $L^1\log(L)^2$. Billard \cite{MR217510} adapted Carleson's arguments to prove almost everywhere convergence of Walsh--Fourier series
for functions in $L^2$.

Fefferman gave an alternative proof of Carleson's theorem \cite{fefferman} via an a priori bound for a certain maximal operator in harmonic analysis, known as the Carleson operator.
Lacey and Thiele pointed out a duality between the approaches by Carleson and Fefferman and presented a more symmetric and self-dual proof \cite{lacey-thiele}.

In recent years, the impact of Carleson's theorem has increased thanks to its connections to areas of mathematics and physics including ergodic theory and scattering theory, and various strengthenings and variants of Fefferman's estimates for Carleson's operator have appeared \cite[\S 1]{CarlesonBlueprint}. For more details about Carleson's theorem and its generalizations, we refer the reader to the surveys by Lacey \cite{MR2091007} and Demeter \cite{MR3334208}.

\subsection{Overview of the formalization project}\label{subsec:overview}
The goal of the Carleson project was to formalize the proof of the metric space Carleson theorem, \Cref{metric-space-Carleson}, as well as the proof of the classical Carleson theorem, \Cref{classical-carleson}, as a corollary of this more general result.

The metric space Carleson theorem is fully proven in \cite{ThieleCarlesonPreprint}, although we give a brief proof sketch in \S\ref{subsec:main_results}.
The original reference \cite{ThieleCarlesonPreprint} was written as a traditional paper for experts in harmonic analysis, but its formalization had from the beginning been envisioned as a large scale collaborative project. This posed a problem, since many potential contributors to the formalization would not have the required expertise to fill in the routine arguments that are usually omitted in a harmonic analysis paper.
To address this challenge, before the formalization started, the authors of \cite{ThieleCarlesonPreprint} wrote a much more detailed document \cite{CarlesonBlueprint} intended as a blueprint for the formalization, which was modified as needed while the formalization process was taking place (see \S\ref{sec:blueprint}).

The Carleson formalization project was publicly announced and opened to contributors in June 2024, and it was completed in July 2025. The core authors of the formalization are the 14 authors of this paper, and smaller contributions were made by the 14 additional people listed in the acknowledgements section. This makes the Carleson project one of the largest scale formalization efforts in mathematics to date.

The formalization was written in Lean~4, making extensive use of Lean's mathematical library \mathlib~\cite{mathlib}.\footnote{\url{https://github.com/leanprover-community/mathlib4}}
The code was developed in a public GitHub repository with an associated website\footnote{\url{https://florisvandoorn.com/carleson/}} that also linked to HTML and PDF versions of the blueprint. Coordination among project collaborators took place using the Lean community Zulip chat.\footnote{\url{https://leanprover.zulipchat.com}} For more details, see \S\ref{sec:org}.

When the project was finished, it comprised around 35.6k lines of Lean code (not including spaces or comments), of which around 9.6k are intended to be upstreamed to the \mathlib library, a process that is still ongoing. Pull requests into \mathlib coming from this project are marked with the \lean{carleson} tag; at the time of writing this article, this includes 111 merged pull requests and 5 more under review.

Throughout the article, we introduce Lean code excerpts from the formalization and use the following clickable symbol\href{https://github.com/fpvandoorn/carleson/blob/AFM-frozen/Carleson/Defs.lean#L298-L300}{\extlink} to link to the corresponding code in our repository or in \mathlib, as appropriate. All code excerpts refer to a fixed, frozen version of this repository\href{https://github.com/fpvandoorn/carleson/tree/AFM-frozen}{\extlink} and the corresponding version of \mathlib (on August 22, 2026).
Some of the code excerpts have been slightly edited for readability.

\subsection{Related work}\label{subsec:rel_work}

Large formalization projects in mathematics have been undertaken in different proof assistants and various areas.
Notable projects formalizing results related to analysis include Immler's verification of the singular hyperbolicity of the Lorenz attractor \cite{ImmlerLorentzSolver} (in Isabelle/HOL),
Immler and Tan's formalization of the Poincaré--Bendixson theorem \cite{ImmlerTanPoincareBendixon} (in Isabelle/HOL),
the Flyspeck project \cite{FlyspeckComplete} formally verifying Hales' proof of the Kepler conjecture (verified using an interplay of results in the Isabelle/HOL and HOL/Light proof assistants), and the formalization of the prime number theorem in HOL Light \cite{HarrisonPNT} (by Harrison), Isabelle/HOL \cite{AvigadPNT} (by Avigad et al), Metamath \cite{CarneiroPNT} (by Carneiro) and Lean (by Kontorovich et al).\footnote{\url{https://alexkontorovich.github.io/PrimeNumberTheoremAnd/}}

Other projects formalizing a single theorem or cluster of theorems involving large collaborations include the formal proof of the Feit--Thompson theorem about finite groups in the Rocq prover \cite{OddOrderTheorem} with a 15 author team, the Liquid Tensor Experiment formalizing an important technical lemma about condensed mathematics in Lean \cites{LTEDefinitions,LTEFinal} (involving 18 authors), the PNT$+$ project formalizing several versions of the prime number theorem using modern tools (with 75 contributors as of September 2026)
and the Equational Theories project classifying certain equational laws on magmas \cite{tao2025equationaltheories} (with 34 authors).
The last project is notable for involving a variety of different tools,
all being used to generate eventual proofs in Lean.
The Flyspeck project was also a large collaborative project, eventually involving 22 authors.

The need for a shared mathematical proof arises naturally in large collaborations.
A common feature is an explicit mathematical \emph{blueprint}, a detailed natural language proof carefully written to be accessible for non-experts. In the Flyspeck project, this proof was written as a book prior to the formalization \cite{HalesFlyspeckBook}. More recent formalization projects (including the Liquid Tensor Experiment and the Sphere Eversion project \cite{SphereEversionCPP}) have used special software (see \S\ref{sec:org}) to maintain such a natural language proof side by side with the formal proof.

Prior to the Carleson project, relatively little harmonic analysis had been formalized in any proof assistant: for example, the Isabelle and HOL Light libraries have an implementation of Fourier series\footnote{\url{https://isa-afp.org/entries/Fourier.html}} (but no more advanced tools such as the real interpolation theorem); the Mizar Mathematical library has a formalization of Minkowski's inequality,\footnote{for infinite sums, \url{https://fm.mizar.org/2005-13/pdf13-1/holder_1.pdf}} but no definition of Fourier series (only discrete Fourier transforms).
Mathlib contains some Fourier analysis thanks to Doll's work on tempered distributions \cite{DollSobolevFourier}.
The Hardy--Littlewood maximal estimate has been formalized in Rocq \cites{AffeldtFTCRocq,AffeldtLebesgueDiffRocq}, but only for functions $\real\to\real$.

\subsection*{Acknowledgements}
The authors are grateful to Asgar Jamneshan, Rajula Srivastava and Christoph Thiele for writing the mathematical paper \cite{ThieleCarlesonPreprint} and the blueprint \cite{CarlesonBlueprint} on which this formalization is based.
We furthermore acknowledge contributions in the form of small formalization additions,
pointing out corrections to the blueprint,
or supplying ideas to the Lean efforts by the following people:
Michel Alexis,
Bolton Bailey,
Julian Berman,
Joachim Breitner,
Martin Dvořák,
Georges Gonthier,
Aaron Hill,
Austin Letson,
Bhavik Mehta,
Eric Paul,
Clara Torres,
Dennis Tsar,
Andrew Yang, and
Ruben van de Velde.

L.B., M.I.d.F.F., L.D., F.v.D., M.R.\ were funded by the Deutsche For\-schungs\-gemein\-schaft (DFG, German Research Foundation) under Germany's Excellence Strategy -- EXC-2047/1 -- 390685813. L.D., F.v.D., M.R.\ were further supported by ERC Synergy Grant 101224275.
J.R.\ was supported in part by NSF grant DMS-2154835 and a HIM fellowship for the Fall 2024 trimester program in Bonn.

\subsection*{Declaration about AI use} Generative AI was used in minor ways in this project, mostly for finding
 results formalized in Lean and for searching the literature. The vast majority of the code is human-written, with the exception of \href{https://github.com/fpvandoorn/carleson/pull/521}{pull request 521}, written after the project was finished, which contains some code generated by the AI agent Aristotle and edited by the PR author. The entirety of this paper is human-written, but AI tools have been used to check for typos and suggest minor fixes.

\section{The metric space Carleson theorem}\label{sec:main_thm}

\subsection{Mathematical background}\label{subsec:background}
The Carleson theorem can be stated for Fourier series or Fourier transforms of a function, but in our formalization we targeted Carleson's theorem for the Fourier series of a periodic function on the real numbers.
Given a $2\pi$-periodic function $f\colon \real \to \complex$, integrable on $[0, 2\pi]$, and $n \in \mathbb{Z}$, define the \emph{$n$-th Fourier coefficient} of $f$ as
\begin{equation*}
    c_n \coloneq \widehat{f}_n \coloneq \frac 1{2\pi}\int_0^{2\pi}f(x) e^{- i nx}\, dx.
\end{equation*}
The Lean version \lean{fourierCoeffOn}\href{https://github.com/leanprover-community/mathlib4/blob/b63493a4746b4651fceb3bcf3a4651cc36a4b8de/Mathlib/Analysis/Fourier/AddCircle.lean#L353-L355}{\extlink} uses \mathlib's formalization of the Bochner integral.

\begin{definition}\label{def:fourier-series}
The \emph{$N$-th partial Fourier sum}\href{https://github.com/fpvandoorn/carleson/blob/AFM-frozen/Carleson/Defs.lean#L43-L44}{\extlink} of $f$ for $x\in\real$ is defined as
\begin{equation*}
    S_N f(x) \coloneq \sum_{n=-N}^N c_n e^{inx}. \
\end{equation*}
The \emph{Fourier series} of $f$ is then given by the formal trigonometric series $\sum_{n=-\infty}^\infty c_n e^{inx}$,
which may or may not converge for a given $x\in\real$.
\end{definition}

The metric space Carleson theorem, \Cref{metric-space-Carleson}, is a so-called restricted weak type estimate for a generalized Carleson operator. One of the generalizations here is that the domain of the functions operated on can be any \emph{doubling metric measure space} $(X,\rho,\mu, a)$.

\begin{definition}\label{def:doubling}
A \emph{doubling metric measure
space}\href{https://github.com/fpvandoorn/carleson/blob/AFM-frozen/Carleson/Defs.lean#L98-L102}{\extlink}
$(X,\rho,\mu, a)$ is a complete and locally compact metric space $(X,\rho)$
equipped with a non-zero locally finite Borel measure $\mu$ that satisfies
the following doubling condition for the natural number $a$. For all $x\in X$ and
all $R>0$ we have
\begin{equation*}
	\mu(B(x,2R))\le 2^a\mu(B(x,R)),
\end{equation*}
where $B(x,R):=\{y\in X: \rho(x,y)<R\}$ denotes the open ball.
\end{definition}

The proof of \Cref{metric-space-Carleson} relies on several tools from real analysis, including the Marcin\-kiewicz real interpolation theorem (see \S\ref{subsubsec:real_interpolation}). Its statement and proof require the following definitions for which we could rely on \mathlib's measure theory (see e.g.\ \cite{vanDoornFormalizedHaarMeasure} for an overview).
Let $E$ be either a normed space over $\mathbb{R}$ endowed with a norm $\|\cdot\|$,
or the extended non-negative real numbers $\ennreal=[0,\infty]$, where $\|\cdot\|$ is defined to be the identity function (see \S\ref{subsec:enorm}). For a function $f$ on a measure space $(X,\mu)$ taking values in $E$, we define its \emph{$L^p$-norm}\href{https://github.com/leanprover-community/mathlib4/blob/b63493a4746b4651fceb3bcf3a4651cc36a4b8de/Mathlib/MeasureTheory/Function/LpSeminorm/Defs.lean#L83-L87}{\extlink} to be
\begin{equation*}
  \|f\|_{L^{p}} \coloneq \big(\int_X\|f(x)\|^p\, d\mu(x)\big) ^ \frac{1}{p}
\end{equation*}
for $0<p<\infty$. For $p=\infty$, the norm $\|f\|_{L^{p}}$ is the essential supremum of $f$, i.e.\ the least value $t$ such that
$\mu\{x\in X \mid t < \|f(x)\| \}=0$. Similarly, we define the \emph{weak $L^p$-norm}\href{https://github.com/fpvandoorn/carleson/blob/AFM-frozen/Carleson/ToMathlib/WeakType.lean#L37-L38}{\extlink} as
\begin{equation*}
  \|f\|_{L^{p,\infty}} \coloneq \sup_{t>0} t \cdot \mu (\{x\in X \mid t < \|f(x)\| \}) ^ \frac{1}{p}
\end{equation*}
for $0<p<\infty$, and $\|f\|_{L^{\infty,\infty}} \coloneq \|f\|_{L^{\infty}}$.
The function $f$ is said to be \emph{in weak $L^p$}\href{https://github.com/fpvandoorn/carleson/blob/AFM-frozen/Carleson/ToMathlib/WeakType.lean#L285-L286}{\extlink} if $f$ is $\mu$-a.e.\ strongly measurable with respect to a fixed topology on $E$ (see e.g.\ \cite{GouezelChangeOfVariables} for motivation of the different notions of measurability in \mathlib) and $\|f\|_{L^{p,\infty}} < \infty$.
An important consequence is that if $f$ is in weak $L^p$, then $\|f(x)\|$ is finite for almost every $x$.\href{https://github.com/fpvandoorn/carleson/blob/AFM-frozen/Carleson/ToMathlib/WeakType.lean#L302}{\extlink}

Given two measure spaces $(X,\mu)$ and $(Y,\nu)$, and two topological spaces $E$ and $F$ endowed with extended norms, we say that an operator $T\colon (X \to E) \to (Y \to F)$ is \emph{of weak type $(p,p')$}\href{https://github.com/fpvandoorn/carleson/blob/AFM-frozen/Carleson/ToMathlib/WeakType.lean#L418-L420}{\extlink} if there exists a constant $C$ such that for any $f\in L^p(X; E)$, we have that $T f$ is $\nu$-a.e.\ strongly measurable and
\begin{equation*}
  \|T f\|_{L^{p',\infty}} \le C \cdot \|f\|_{L^p}.
\end{equation*}
Similarly, the operator $T$ is of \emph{strong type $(p, p')$}\href{https://github.com/fpvandoorn/carleson/blob/AFM-frozen/Carleson/ToMathlib/WeakType.lean#L432-L435}{\extlink} if there exists a constant $C$ such that for any $f\in L^p(X; E)$, we have that $T f$ is $\nu$-a.e.\ strongly measurable and
\begin{equation*}
  \|T f\|_{L^{p'}} \le C \cdot \|f\|_{L^p}.
\end{equation*}

\subsection{Statement of the main results}\label{subsec:main_results}

\subsubsection{Classical Carleson}\label{subsubsec:classical-result}

The namesake of the formalization project is Carleson's classical result on the pointwise convergence of Fourier series.
The specific statement we formalized is as follows. Note that in our formalization we require that the given function is continuous.

\begin{theorem}[classical Carleson \cite{carleson}\href{https://github.com/fpvandoorn/carleson/blob/AFM-frozen/Carleson/Classical/ClassicalCarleson.lean\#L131-L132}{\extlink}]
    \label{classical-carleson}
    Let $f$ be a $2\pi$-periodic complex-valued continuous function on $\mathbb{R}$.
    Then for almost all $x \in \mathbb{R}$, the limit
    $\lim_{N\to\infty} S_N f(x)$ exists and is equal to $f(x)$.
\end{theorem}

The metric space Carleson theorem is more complicated to state, but gives bounds for certain \emph{generalized Carleson operators}.
Before we state this theorem, we state an intermediate result
bounding a specific Carleson operator for functions on $\real$.
Like most proofs of Carleson's theorem, we use this bound to prove \Cref{classical-carleson}.
The specific version of the bound that we formalized is as follows.

\begin{theorem}[real Carleson\href{https://github.com/fpvandoorn/carleson/blob/AFM-frozen/Carleson/Classical/CarlesonOnTheRealLineContinuous.lean\#L29-L31}{\extlink}]\label{real-Carleson}
    Let $F,G$ be Borel subsets of $\real$ with finite measure. Let $f$ be a bounded measurable function on $\real$ with $|f|\le \mathbf{1}_F$. Then
\begin{equation*}
    \int _G Tf(x) \, dx \le 2^{474 \cdot 4^3} |F|^{\frac 12} |G|^{\frac 12},
\end{equation*}
where
\begin{equation*}
    T f(x)\coloneq\sup_{n\in \mathbb{Z}}
    \sup_{r>0}\left|\int_{r<|x-y|<1} f(y)\kappa(x-y) e^{iny}\, dy\right|,
\end{equation*}
and the \emph{Hilbert kernel} $\kappa\colon\real\to\complex$ is defined as
\begin{equation*}
    \kappa(x)\coloneq\begin{cases}
        0 & \text{if $x = 0$ or $|x| \ge 1$} \\
        \frac { 1-|x|}{1-e^{ix}} & \text{if $0 < |x| < 1$}.
    \end{cases}
\end{equation*}
\end{theorem}

\subsubsection{Preliminary definitions}\label{subsubsec:definitions}

In this section we state the definitions used in the metric space Carleson theorem.
Let $a \ge 4$ be a natural number.
\begin{definition}
A collection $\Mf$ of real-valued continuous functions on the doubling metric measure space $(X,\rho,\mu,a)$ is called \emph{compatible}\href{https://github.com/fpvandoorn/carleson/blob/AFM-frozen/Carleson/Defs.lean#L167-L185}{\extlink}
if there is a point $o\in X$ where all the functions are equal to $0$,
and if there exists for each open ball $B \subseteq X$ a metric $d_B$ on $\Mf$,
such that the following five properties are satisfied.
\begin{itemize}\label{def:compatible}
\item For every open ball $B \subseteq X$ and any $\mfa, \mfb \in \Mf$,
\begin{equation}\label{osccontrol}
    \sup_{x,y\in B}|\mfa(x)-{\mfa(y)}- \mfb(x)+{\mfb(y)}| \le d_{B}(\mfa,\mfb).
\end{equation}
\item For any two open balls $B_1=B(x_1,R)$, $B_2= B(x_2,2R)$ in $X$ with $x_1\in B_2$ and any $\mfa,\mfb\in \Mf$,
\begin{equation}\label{firstdb}
    d_{B_2}(\mfa,\mfb) \le 2^a d_{B_1}(\mfa,\mfb).
\end{equation}
\item For any two open balls $B_1, B_2$ in $X$ with $B_1 \subseteq B_2$ and any $\mfa, \mfb \in \Mf$
\begin{equation}\label{monotonedb}
    d_{B_1}(\mfa,\mfb) \le d_{B_2}(\mfa, \mfb).
\end{equation}
\item For any two open balls
$B_1=B(x_1,R)$, $B_2= B(x_2,2^aR)$
with $B_1\subseteq B_2$, and $\mfa,\mfb\in \Mf$,
\begin{equation}\label{seconddb}
    2d_{B_1}(\mfa,\mfb) \le d_{B_2}(\mfa,\mfb).
\end{equation}
\item For every open ball $B$ in $X$ and every open $d_B$-ball $\tilde B$ of radius $2R$ in $\Mf$, there is a collection $\mathcal{B}$ of
at most $2^a$ many open $d_B$-balls of radius $R$ covering $\tilde B$, that is,
\begin{equation}\label{thirddb}
    \tilde B\subseteq \bigcup \mathcal{B}.
\end{equation}
\end{itemize}
Elements of $\Theta$ are called modulation functions.
\end{definition}

\begin{definition}
A compatible collection $\Mf$ is called \emph{cancellative}\href{https://github.com/fpvandoorn/carleson/blob/AFM-frozen/Carleson/Defs.lean#L191-L197}{\extlink} if
for any open ball $B$ in $X$ of radius $R$, any Lipschitz function $\varphi\colon X\to \complex$
supported on $B$, and any $\mfa,\mfb\in \Mf$ we have
\begin{equation*}
	\Big|\int_B e(\mfa(x)-{\mfb(x)}) \varphi(x) d\mu(x)\Big|\le 2^a \mu(B)\|\varphi\|_{\Lip(B)}
	(1+d_B(\mfa,\mfb))^{-\frac{1}{a}},
\end{equation*}
where $e(r) = e^{ir}$ and $\|\cdot\|_{\Lip(B)}$ denotes the inhomogeneous Lipschitz norm on $B$:
$$
\|\varphi\|_{\Lip(B)} = \sup_{x \in B} |\varphi(x)| + R \sup_{x,y \in B, x \neq y} \frac{|\varphi(x) - \varphi(y)|}{\rho(x,y)}.
$$
\end{definition}

\begin{definition}\label{def:one-sided-kernel}
A \emph{one-sided Calderón--Zygmund kernel}\href{https://github.com/fpvandoorn/carleson/blob/AFM-frozen/Carleson/Defs.lean#L250-L254}{\extlink} $K$ on the doubling metric measure space $(X, \rho, \mu, a)$
is a measurable function $K:X\times X\to \mathbb{C}$ such that for all $x,y,y'\in X$ with $x\neq y$, we have
\begin{equation*}
	|K(x,y)| \leq \frac{2^{a^3}}{V(x,y)}
\end{equation*}
and if $2\rho(y,y') \leq \rho(x,y)$, then
\begin{equation*}
	|K(x,y) - K(x,y')| \leq \left(\frac{\rho(y,y')}{\rho(x,y)}\right)^{\frac{1}{a}}\frac{2^{a^3}}{V(x,y)},
\end{equation*}
where $V(x,y):=\mu(B(x,\rho(x,y)))$.
\end{definition}

\begin{definition}\label{def:two-sided-kernel}
A \emph{two-sided Calderón--Zygmund kernel}\href{https://github.com/fpvandoorn/carleson/blob/AFM-frozen/Carleson/TwoSidedCarleson/WeakCalderonZygmund.lean#L17-L20}{\extlink} is a one-sided Calderón--Zygmund kernel $K$ such that $(x,y)\mapsto K(y,x)$ is also a one-sided Calderón--Zygmund kernel. Explicitly, this means that for all $x,x',y\in X$ with $x\neq y$ and $2\rho(x,x') \leq \rho(x,y)$,
\begin{equation*}
	|K(x,y) - K(x',y)| \leq \left(\frac{\rho(x,x')}{\rho(x,y)}\right)^{\frac{1}{a}}\frac{2^{a^3}}{V(x,y)}.
\end{equation*}
\end{definition}

Let $f\colon X \to \complex$ be a bounded, measurable function supported on a set of finite measure.
Define the \emph{maximally truncated non-tangential singular integral}\href{https://github.com/fpvandoorn/carleson/blob/AFM-frozen/Carleson/Defs.lean#L211-L212}{\extlink} $T_{*}$ associated with $K$ by
\begin{equation}
    \label{def-tang-unm-op}
    T_{*}f(x):=\sup_{0 < R_1 < R_2} \sup_{\rho(x,x')<R_1} \left|\int_{R_1< \rho(x',y) < R_2} K(x',y) f(y) \, \mathrm{d}\mu(y) \right|.
\end{equation}
Additionally, define for $r > 0$, $x\in X$ the \emph{Calderon--Zygmund operator}\href{https://github.com/fpvandoorn/carleson/blob/AFM-frozen/Carleson/TwoSidedCarleson/Basic.lean#L15-L16}{\extlink} $T_r$ as
\begin{equation}
\label{def-T-r}
T_r f(x):= \int_{r\le\rho(x,y)} K(x,y) f(y) \, d\mu(y) = \int_{X\setminus B(x,r)} K(x,y) f(y) \, d\mu(y).
\end{equation}
Finally, define the \emph{generalized Carleson operator}\href{https://github.com/fpvandoorn/carleson/blob/AFM-frozen/Carleson/Defs.lean#L228-L229}{\extlink} $T$ by
\begin{equation}
    \label{def-main-op}
    Tf(x):=\sup_{\mfa\in\Mf} \sup_{0 < R_1 < R_2}\left| \int_{R_1 < \rho(x,y) < R_2} K(x,y) f(y) e(\mfa(y)) \, \mathrm{d}\mu(y) \right|,
\end{equation}
where $e(r)=e^{ir}$.

For a Borel function $\tQ:X\to \Mf$, and $\mfa \in \Mf$ and $x\in X$ define
\begin{equation*}
    R_{\tQ}(\mfa,x)=\sup\{r>0:d_{B(x,r)}(\mfa, \tQ(x))<1\}
\end{equation*}
and define further the \emph{linearized maximally truncated nontangential Calderon–Zygmund operator}\href{https://github.com/fpvandoorn/carleson/blob/38587d19c6c4c38a179b68476f717cad2706c51e/Carleson/Defs.lean#L204-L208}{\extlink} $T_{\tQ}^\mfa$ as
\begin{equation}\label{def-lin-star-op}
    T_{\tQ}^\mfa f(x):=\sup_{0<R_1<R_2} \ \sup_{\rho(x,x')<R_1}
    \left|\int_{R_1< \rho(x',y) < \min\{R_2, R_{\tQ}(\mfa,x')\}} K(x',y) f(y) \, \mathrm{d}\mu(y) \right|.
\end{equation}
Finally, define the \emph{linearized generalized Carleson operator}\href{https://github.com/fpvandoorn/carleson/blob/38587d19c6c4c38a179b68476f717cad2706c51e/Carleson/Defs.lean#L222-L224}{\extlink} $T_\tQ$ by
\begin{equation}\label{def-lin-main-op}
    T_\tQ f(x):= \sup_{0 < R_1 < R_2}\left| \int_{R_1 < \rho(x,y) < R_2} K(x,y) f(y) e(\tQ(x)(y)) \, \mathrm{d}\mu(y) \right|,
\end{equation}
where again $e(r)=e^{ir}$.

\subsubsection{Metric space Carleson}\label{subsubsec:metric-carleson}

We prove three versions of the main result.
The first version is the following restricted weak type estimate for $T$ in the range $1<q\le 2$,
which by interpolation techniques recovers $L^q$ estimates for the open range $1<q<2$.
\begin{theorem}[metric space Carleson\href{https://github.com/fpvandoorn/carleson/blob/AFM-frozen/Carleson/MetricCarleson/Main.lean\#L203-L208}{\extlink}]
\label{metric-space-Carleson}
    For all integers $a \ge 4$ and real numbers $1<q\le 2$ the following holds.
    Let $(X,\rho,\mu,a)$ be a doubling metric measure space.
    Let $\Mf$ be a cancellative compatible collection of functions and let $K$ be a one-sided Calderón--Zygmund kernel on $(X,\rho,\mu,a)$.
    Assume that for every bounded measurable function $g$ on $X$ supported on a set of finite measure we have
    \begin{equation}\label{nontanbound}
        \|T_{*}g\|_{2} \leq 2^{a^3} \|g\|_2,
    \end{equation}
    where $T_{*}$ is defined in \eqref{def-tang-unm-op}.
    Then for all Borel subsets $F$ and $G$ of $X$ and all Borel functions $f\colon X\to \complex$ with
    $|f|\le \mathbf{1}_F$, we have, with $T$ defined in \eqref{def-main-op},
    \begin{equation*}
        \int_{G} T f \, \mathrm{d}\mu \leq \frac{2^{443a^3}}{(q-1)^6} \mu(G)^{1-\frac{1}{q}} \mu(F)^{\frac{1}{q}}.
    \end{equation*}
\end{theorem}

In some applications, such as the Walsh-case of Carleson's theorem \cite{MR217510}, the kernel $K$ naturally depends also on the modulation functions $\mfa$.
The fact that we don't assume Hölder continuity of the kernel $K$ in the first argument allows us to also capture this situation.

The second version of our main result has a slightly weaker replacement of the assumption \eqref{nontanbound}.
In fact, it implies the Walsh-case of Carleson's theorem.

\begin{theorem}[linearized metric space Carleson\href{https://github.com/fpvandoorn/carleson/blob/AFM-frozen/Carleson/MetricCarleson/Linearized.lean\#L109-L115}{\extlink}]
\label{linearized-metric-Carleson}
    For all integers $a \ge 4$ and real numbers $1<q\le 2$ the following holds.
    Let $(X,\rho,\mu,a)$ be a doubling metric measure space. Let $\Mf$ be a
    cancellative compatible collection of functions.
    Let $\tQ:X\to \Mf$ be a Borel function with finite range.
    Let $K$ be a one-sided Calderón--Zygmund kernel on $(X,\rho,\mu,a)$. Assume that for every $\mfa\in \Mf$ and every bounded measurable function $g$ on $X$ supported on a set of finite measure we have
    \begin{equation}\label{linnontanbound}
        \|T_{\tQ}^\mfa g\|_{2} \leq 2^{a^3} \|g\|_2,
    \end{equation}
    where $T_{\tQ}^\mfa$ is defined in \eqref{def-lin-star-op}.
    Then for all Borel subsets $F$ and $G$ of $X$ and all Borel functions $f\colon X\to \complex$ with
    $|f|\le \mathbf{1}_F$, we have, with $T_\tQ$ defined in \eqref{def-lin-main-op},
    \begin{equation*}
        \int_{G} T_\tQ f \, \mathrm{d}\mu \le \frac{2^{443a^3}}{(q-1)^6} \mu(G)^{1-\frac{1}{q}} \mu(F)^{\frac{1}{q}}.
    \end{equation*}
\end{theorem}

The third version of the main result involves two-sided Calderón--Zygmund kernels, whose additional regularity allows us to weaken one of the assumptions to the family of operators $T_r$, which is easier to work with in applications.
In particular, the two-sided case can be applied to derive the classical Carleson theorem.

\begin{theorem}[two-sided metric space Carleson\href{https://github.com/fpvandoorn/carleson/blob/AFM-frozen/Carleson/TwoSidedCarleson/MainTheorem.lean\#L30-L35}{\extlink}]
    \label{two-sided-metric-space-Carleson}
        For all  integers $a \ge  4$ and real numbers $1<q\le 2$ the following holds.
        Let $(X,\rho,\mu,a)$ be a doubling metric measure space. Let $\Mf$ be a
        cancellative compatible collection of functions and let $K$ be a two-sided Calderón--Zygmund kernel on $(X,\rho,\mu,a)$. Assume that for every bounded measurable function $g$ on $X$ supported on a set of finite measure and all $r>0$ we have
      \begin{equation*}
            \|T_r g\|_{2} \leq 2^{a^3} \|g\|_2.
        \end{equation*}
        Then for all Borel subsets $F$ and $G$ of $X$ and
        all Borel functions $f\colon X\to \complex$ with
        $|f|\le \mathbf{1}_F$, we have, with $T$ defined in \eqref{def-main-op},
      \begin{equation*}
            \int_{G} T f \, \mathrm{d}\mu \leq \frac{2^{474a^3}}{(q-1)^6} \mu(G)^{1-\frac{1}{q}} \mu(F)^{\frac{1}{q}}.
        \end{equation*}
\end{theorem}

\subsection{Proof roadmap}\label{subsec:proof-sketch}

In order to highlight choices made in the proof and the formalization it is practical to give a brief roadmap
of the proof of the main results stated above.
Of these three, the linearized version \Cref{linearized-metric-Carleson} implies the others. It is proved by first reducing it to a finitary version, which is then reduced to a discrete version.

Given the variables and assumptions of \Cref{linearized-metric-Carleson}, define $D := 2^{100a^2}$ and $\kappa := 2^{-10a}.$
Let $\psi\colon\real \to\real$ be the unique compactly supported, piecewise linear, continuous function with corners precisely at $\frac 1{4D}$, $\frac 1{2D}$, $\frac 14$ and $\frac 12$ which satisfies
\begin{equation*}
    \sum_{s\in \mathbb{Z}} \psi(D^{-s}x)=1
\end{equation*}
for all $x>0$. Note that the above sum has at most two non-zero terms, and that this function vanishes outside $[\frac1{4D},\frac 12]$, is constant one on
$[\frac1{2D},\frac 14]$, and is Lipschitz with constant $4D$. The finitary Carleson theorem is a variant of \Cref{linearized-metric-Carleson} where the integration domain and the support and range of $f$ are assumed to be finite. \Cref{linearized-metric-Carleson} follows from it by standard limiting arguments.

\begin{theorem}[finitary Carleson\href{https://github.com/fpvandoorn/carleson/blob/AFM-frozen/Carleson/FinitaryCarleson.lean\#L110-L113}{\extlink}]
\label{finitary-Carleson}
Let ${\sigma_1},\sigma_2\colon X\to \mathbb{Z}$ be measurable functions with finite range and ${\sigma_1}\leq \sigma_2$. Let $F,G$ be bounded Borel subsets of $X$. Then there is a Borel subset $G'$ of $X$ with $2\mu(G')\leq \mu(G)$ such that
for all Borel functions $f\colon X\to \complex$ with $|f|\le \mathbf{1}_F$ we have
\begin{equation*}
    \int_{G \setminus G'} \left|\sum_{s={\sigma_1}(x)}^{{\sigma_2}(x)} \int K_s(x,y) f(y) e(\tQ(x)(y)) \, \mathrm{d}\mu(y) \right| \mathrm{d}\mu(x)
    \le \frac{2^{442a^3}}{(q-1)^5} \mu(G)^{1-\frac{1}{q}} \mu(F)^{\frac 1 q},
\end{equation*}
where for $s \in \mathbb{Z}, K_s(x,y) := K(x,y)\psi(D^{-s}\rho(x,y))$.
\end{theorem}

To reduce this finitary version of the Carleson theorem to the discrete version, we make divisions in both the spatial and frequency domains. The spatial discretization will be achieved by a \emph{grid structure}, after which a subordinate frequency discretization is encoded by a \emph{tile structure}.

\begin{definition}\label{def:grid-structure}
A \emph{dyadic cube} is a pair $(I, k)$ of a Borel subset $I$ of $X$ and an integer scale $k \in [-S, S]$ for some fixed natural number $S$. Usually we abuse notation and write $I$ for a dyadic cube $(I, k)$. We write $I \le J$ to mean that for two cubes $(I, k)$ and $(J, l)$, we have $I \subseteq J$ and $k \le l$.

A \emph{grid structure}\href{https://github.com/fpvandoorn/carleson/blob/AFM-frozen/Carleson/GridStructure.lean#L20-L44}{\extlink} on $X$ is a triple $(\mathcal{D}, c, s)$ consisting of a finite set $\mathcal{D}$ of dyadic cubes, a function $c\colon \mathcal{D} \to X$ called the \emph{center function}, and a surjective function $s\colon \mathcal{D} \to [-S, S]$, $(I, k) \mapsto k$, called the \emph{scale function}, that together satisfy the conditions \ref{it:coverdyadic}--\ref{it:eq-small-boundary}:
\begin{enumerate}[label=(\roman*)]
\item \label{it:coverdyadic} For every dyadic cube $I$ and every $-S \leq k < s(I)$,
\begin{equation*}
    I\subseteq \bigcup_{J\in \mathcal {D}: s(J)=k}J.
\end{equation*}
\item \label{it:dyadicproperty} For all non-disjoint dyadic cubes $I, J$ with $s(I) \leq s(J)$, we have that
\begin{equation*}
    I\le J.
\end{equation*}
\item \label{it:subsetmaxcube} There exists a cube $I_0 \in \mathcal{D}$ with scale $s(I_0) = S$ and center $c(I_0) = o$ (the point where all $\mfa \in \Mf$ vanish), and for all dyadic cubes $J$ it holds that $J \le I_0$.
\item \label{it:eq-vol-sp-cube} Every dyadic cube $I$ fits a ball with center $c(I)$ and is approximately of size $D^{s(I)}$ in the sense that
\begin{equation}
    \label{eq-vol-sp-cube}
    B(c(I), \frac{1}{4} D^{s(I)}) \subseteq I \subseteq B(c(I), 4 D^{s(I)}).
\end{equation}
\item \label{it:eq-small-boundary} For every dyadic cube $I$ and every $t$ with $t D^{s(I)} \geq D^{-S}$, the boundary of $I$ is small in the sense that
\begin{equation*}
    \mu(\{x \in I \mid \rho(x, X \setminus I) \leq t D^{s(I)}\}) \le 2 t^\kappa \mu(I).
\end{equation*}
\end{enumerate}
\end{definition}

We now further subdivide in the frequency domain $\Theta$, intuitively partitioning for every dyadic cube the collection of functions $\Theta$ into pieces called tiles. More precisely, given a grid structure $(\mathcal{D}, c, s)$, a \emph{tile structure}\href{https://github.com/fpvandoorn/carleson/blob/AFM-frozen/Carleson/TileStructure.lean#L54-L65}{\extlink} for this grid structure is a tuple $(\fP,\scI,\fc,\fcc,\pc,\ps)$ where $\fP$ is a finite set with elements called tiles. The function $\fc\colon \fP \to \mathcal{P}(\Mf)$ maps every tile to a subset of the collection of functions $\Theta$, and the function $\fcc\colon \fP \to \tQ(X)$ plays the role of center function for the selection of functions. The association of dyadic cubes to a collection of tiles is encoded in the reverse direction by a surjective function $\scI: \fP \to \mathcal{D}$ that merely associates to every tile a dyadic cube and the functions $\pc\colon \fP \to X$ and $\ps\colon \fP \to [-S, S]$ that just yield the center and scale of the given cube. The above functions need to satisfy the following conditions.

\begin{enumerate}[label=\roman*.]
\item \label{it:omega-disjoint} For every dyadic cube $I$ and every pair of tiles $\fp \neq \fq$ in $\scI^{-1}(\{I\})$, the sets $\fc(\fp)$ and $\fc(\fq)$ are disjoint.
\item For every dyadic cube $I$, we have that
\begin{equation*}
    \tQ(X)\subseteq \bigcup_{\fp\in \scI^{-1}(\{I\})}\fc(\fp).
\end{equation*}
\item For all tiles $\fp, \fq$ with $\scI(\fp) \le \scI(\fq)$ it holds that
\begin{equation*}
    \fc(\fq)\subseteq \fc(\fp) \quad\text{or}\quad \fc(\fq)\cap \fc(\fp) = \emptyset.
\end{equation*}
\item For every tile $\fp$, the set of functions $\fc(\fp)$ contains a ball with center $\fcc(\fp)$ and is of size approximately $D^{\ps(\fp)}$ in the sense that
\begin{equation*}
    B_{\fp}(\fcc(\fp), 0.2) \subseteq \fc(\fp) \subseteq B_{\fp}(\fcc(\fp),1),
\end{equation*}
where
\begin{equation*}
    B_{\fp} (\mfa, R) := \{\mfb \in \Mf \mid d_{\fp}(\mfa, \mfb) < R\},
\end{equation*}
with
\begin{equation*}
    d_{\fp} := d_{B(\pc(\fp),\frac 14 D^{\ps(\fp)})}.
\end{equation*}
\end{enumerate}

With the above definitions, it can be shown that the finitary Carleson theorem reduces to the discrete Carleson theorem.

\begin{theorem}[discrete Carleson\href{https://github.com/fpvandoorn/carleson/blob/AFM-frozen/Carleson/Discrete/MainTheorem.lean\#L44-L48}{\extlink}]
\label{discrete-Carleson}
Let $(\mathcal{D}, c, s)$ be a grid structure and $(\fP,\scI,\fc,\fcc,\pc,\ps)$
a tile structure for this grid structure.
Define for every tile $\fp\in \fP$
\begin{equation*}
    E(\fp)=\{x\in \scI(\fp): \tQ(x)\in \fc(\fp) , {\sigma_1}(x)\le \ps(\fp)\le {\sigma_2}(x)\}
\end{equation*}
and
\begin{equation*}
    T_{\fp} f(x)= \mathbf{1}_{E(\fp)}(x) \int K_{\ps(\fp)}(x,y) f(y) e(\tQ(x)(y)-\tQ(x)(x))\, d\mu(y).
\end{equation*}
Then there exists a Borel subset $G'$ of $X$ with $2\mu(G') \leq \mu(G)$ such that for all Borel functions $f\colon X\to \complex$ with $|f|\le \mathbf{1}_F$
we have
\begin{equation*}
    \int_{G \setminus G'} \left| \sum_{\fp \in \fP} T_{\fp} f (x) \right| \, \mathrm{d}\mu(x) \le \frac{2^{442a^3}}{(q-1)^5} \mu(G)^{1-\frac{1}{q}} \mu(F)^{\frac{1}{q}}.
\end{equation*}
\end{theorem}

The proof of this last result is the heart of the matter, and is very technical.
The full proof can be found in sections 5, 6 and 7 of the blueprint \cite{CarlesonBlueprint};
let us give an idea of the work required.

Firstly, cubes are classified by what fraction of them is covered by the set $G$.
This means that we partition the set of cubes into collections $\mathcal{C}(k)$
for which this fraction is in the interval $(2^{-(k+1)},2^{-k}]$.
This classification lifts via the function $\scI$ to the tiles.
Secondly, we make a subclassification of tiles $\fp$
by the fraction of the points $x$ in $\scI(\fp)$ for which $d_{\fp}(\tQ(x),\fcc(\fp))$ is small.
As a further refinement, some tiles are \emph{marked},
and we classify tiles by how many marked tiles are ``near'' a given tile.
Let us note that these characterizations also take into account the orders on cubes and tiles.
The order on tiles is given by $\fp\le \fp'$ iff
$\scI(\fp)\le \scI(\fp')$ and $\Omega(\fp)\supseteq \Omega(\fp')$.

In each class in the final classification, we mark certain tiles as so-called \emph{tree tops}.
Then we eliminate five separate subclasses of tiles from the class.
Four of these subclasses are antichains with respect to the aforementioned order on the tiles.
The main result of Section~6 of the blueprint gives good bounds for the
Carleson operator on an antichain $\mathfrak{A}\subseteq\fP$, i.e.\ for the operator
$$\sum_{\fp\in\mathfrak{A}}T_{\fp}.$$
The fifth subclass is small, and we can choose $G'$
so that it contains $\scI(\fp)$ for all tiles $\fp$ in this subclass.

The remaining tiles can be partitioned into a small number of so-called \emph{forests}\href{https://github.com/fpvandoorn/carleson/blob/AFM-frozen/Carleson/Forest.lean#L21-L33}{\extlink}
(not in the graph-theoretic sense)
with trees that each contain one of the aforementioned tree tops.
The main result of Section~7 of the blueprint estimates the Carleson operator on such a forest.
Section~5 gives the characterization of tiles,
and puts these results together to prove the discrete Carleson theorem.

\subsection{Stating the main results in Lean}\label{subsec:lean-statements}

Starting from the material already available in Lean's mathematical library \mathlib, it is easy to state the classical Carleson theorem in Lean.
\begin{lstlisting}[mathescape]
theorem classical_carleson {f : ℝ → ℂ}
  (cont_f : Continuous f) (periodic_f : f.Periodic (2 * π)) :
  ∀$^m$ x, Tendsto (partialFourierSum · f x) atTop (nhds (f x))
\end{lstlisting}

The statement of the hypotheses ($f \colon \real \to \complex$ is a continuous $2\pi$-periodic function) closely resembles natural language, but understanding the conclusion requires some familiarity with Lean and \mathlib. There, $\forall^m x$ means ``for every real number $x$ not in a set of measure zero'', and
\lean{Tendsto (partialFourierSum · f x) atTop (nhds (f x))} is Lean's idiomatic way to express that the sequence $S_N f(x)$ of $N$-th partial Fourier sums converges to $f(x)$ when $N$ tends to infinity (see \cite[\S11.1]{mathematics_in_lean} for more information on \mathlib's use of filters to represent limits).

The definition \lean{partialFourierSum} encodes the $N$-th partial Fourier sum of $f\colon\real \to\complex$. Here, \lean{Finset.Icc}\href{https://github.com/leanprover-community/mathlib4/blob/b63493a4746b4651fceb3bcf3a4651cc36a4b8de/Mathlib/Order/Interval/Finset/Defs.lean#L280}{\extlink}\footnote{In any preorder $\alpha$, one can define intervals which on each side can be either open (o), closed (c) or infinite (i). For instance, \lean{Set.Ico a b}\href{https://github.com/leanprover-community/mathlib4/blob/b63493a4746b4651fceb3bcf3a4651cc36a4b8de/Mathlib/Order/Interval/Set/Defs.lean\#L59}{\extlink} denotes the interval $[a, b)$ and \lean{Set.Iic b}\href{https://github.com/leanprover-community/mathlib4/blob/b63493a4746b4651fceb3bcf3a4651cc36a4b8de/Mathlib/Order/Interval/Set/Defs.lean\#L44}{\extlink} denotes $(-\infty, b]$. If moreover $\alpha$ is a \lean{LocallyFiniteOrder}\href{https://github.com/leanprover-community/mathlib4/blob/b63493a4746b4651fceb3bcf3a4651cc36a4b8de/Mathlib/Order/Interval/Finset/Defs.lean\#L108-L124}{\extlink} (i.e., all bounded intervals in $\alpha$ are finite), \lean{Finset.Icc, Finset.Ico, Finset.Ioc} and \lean{Finset.Ioo} provide \lean{Finset}\href{https://github.com/leanprover-community/mathlib4/blob/b63493a4746b4651fceb3bcf3a4651cc36a4b8de/Mathlib/Data/Finset/Defs.lean\#L72-L79}{\extlink} versions of these intervals.
} denotes the finite set of integers $\{-N, -N + 1, \cdots, N\}$, \lean{fourierCoeffOn Real.two_pi_pos f n}\href{https://github.com/leanprover-community/mathlib4/blob/b63493a4746b4651fceb3bcf3a4651cc36a4b8de/Mathlib/Analysis/Fourier/AddCircle.lean#L351-L355}{\extlink} is the
$n$-th Fourier coefficient of $f$ on $[0, 2\pi]$, and \lean{fourier n (x : AddCircle (2 * π))}\href{https://github.com/leanprover-community/mathlib4/blob/b63493a4746b4651fceb3bcf3a4651cc36a4b8de/Mathlib/Analysis/Fourier/AddCircle.lean#L123-L125}{\extlink} denotes $e^{inx}$.
\begin{lstlisting}[mathescape]
def partialFourierSum (N : ℕ) (f : ℝ → ℂ) (x : ℝ) : ℂ :=
  ∑ n ∈ Finset.Icc (-(N : ℤ)) N,
    fourierCoeffOn Real.two_pi_pos f n * fourier n (x : AddCircle (2 * π))
\end{lstlisting}

In contrast to the classical result, to be able to formalize the statements of the metric and linearized versions of Carleson's theorem, we first needed to formalize a good amount of new concepts.
The first such concept is a doubling metric measure space (see Definition~\ref{def:doubling}), which we denote by \lean{MeasureTheory.DoublingMeasure}\href{https://github.com/fpvandoorn/carleson/blob/AFM-frozen/Carleson/Defs.lean#L98-L102}{\extlink}. The constant $A$ appearing in this definition will typically be $2^a$.
\begin{lstlisting}[mathescape]
class MeasureTheory.Measure.IsDoubling {X : Type*} [MeasurableSpace X]
    [PseudoMetricSpace X] (μ : Measure X) (A : outParam ℝ≥0) : Prop where
  measure_ball_two_le_same : ∀ (x : X) r, μ (ball x (2 * r)) ≤ A * μ (ball x r)

class MeasureTheory.DoublingMeasure (X : Type*) (A : outParam ℝ≥0) [PseudoMetricSpace X] extends
  CompleteSpace X, LocallyCompactSpace X, MeasureSpace X, BorelSpace X,
  IsLocallyFiniteMeasure (volume : Measure X),
  IsDoubling (volume : Measure X) A, NeZero (volume : Measure X) where
\end{lstlisting}

Next, we need to introduce the family $\Mf$ of compatible functions satisfying the properties \eqref{osccontrol}, \eqref{firstdb}, \eqref{monotonedb}, \eqref{seconddb}, and \eqref{thirddb}. We do this in two steps, by first creating a class
\lean{FunctionDistances}\href{https://github.com/fpvandoorn/carleson/blob/AFM-frozen/Carleson/Defs.lean#L113-L121}{\extlink} which bundles a family of continuous functions and then endowing these functions with the five properties in the class \lean{CompatibleFunctions}.\href{https://github.com/fpvandoorn/carleson/blob/AFM-frozen/Carleson/Defs.lean#L167-L185}{\extlink} The assumption that the family is cancellative is added as a \lean{Prop}-valued class \lean{IsCancellative}.\href{https://github.com/fpvandoorn/carleson/blob/AFM-frozen/Carleson/Defs.lean#L191-L197}{\extlink}

Something to note here is that we created type synonyms to be able to endow $\Mf$ with different metric space structures, each coming from a ball in $X$, and we provided the custom notation \lean|dist_{x, r}| for the distance corresponding to the ball $B(x, r)$.

\begin{lstlisting}[mathescape]
variable {k X : Type*} {A : ℕ} [RCLike k] [PseudoMetricSpace X]

def WithFunctionDistance (x : X) (r : ℝ) := Θ X

instance {x : X} {r : ℝ} [d : FunctionDistances k X] :
  PseudoMetricSpace (WithFunctionDistance x r) := d.metric x r

notation3 "dist_{" x ", " r "}" => @dist (WithFunctionDistance x r) _
\end{lstlisting}

With this notation, the Lean statement of the compatibility properties is quite close to the informal one. For instance, equation \eqref{monotonedb} is written as
\begin{lstlisting}[mathescape]
lemma cdist_mono {x₁ x₂ : X} {r₁ r₂ : ℝ} {f g : Θ X}
  (h : ball x₁ r₁ ⊆ ball x₂ r₂) : dist_{x₁, r₁} f g ≤ dist_{x₂, r₂} f g
\end{lstlisting}

The class \lean{IsOneSidedKernel} captures the notion of one-sided Calderón--Zygmund kernel (Definition~\ref{def:one-sided-kernel}).

\begin{lstlisting}[mathescape]
class IsOneSidedKernel (a : outParam ℕ) (K : X → X → ℂ) : Prop where
  measurable_K : Measurable (uncurry K)
  norm_K_le_vol_inv (x y : X) : ‖K x y‖ ≤ C_K a / Real.vol x y
  norm_K_sub_le {x y y' : X} (h : 2 * dist y y' ≤ dist x y) : ‖K x y - K x y'‖ ≤
    (dist y y' / dist x y) ^ (a : ℝ)⁻¹ * (C_K a / Real.vol x y)
\end{lstlisting}

To avoid constant repetition of hypotheses, we create a class \lean{KernelProofData}\href{https://github.com/fpvandoorn/carleson/blob/AFM-frozen/Carleson/Defs.lean#L261-L266}{\extlink} which contains data used in most of sections 2 through 7 of the blueprint: a doubling measure on $X$, a natural number greater than three, a collection of compatible functions, and a one-sided Calderón--Zygmund kernel. See \S\ref{subsec:proof_data} for a more detailed discussion of this design choice.
\begin{lstlisting}[mathescape]
class KernelProofData {X : Type*} (a : outParam ℕ) (K : outParam (X → X → ℂ))
    [PseudoMetricSpace X] where
  d : DoublingMeasure X (defaultA a)
  four_le_a : 4 ≤ a
  cf : CompatibleFunctions ℝ X (defaultA a)
  hcz : IsOneSidedKernel a K
\end{lstlisting}

The statement of the metric space Carleson theorem also requires fixing two positive real numbers $q, q'$ satisfying the equality $q^{-1} + q'^{-1} = 1$, encoded by definition \lean{HolderConjugate}\href{https://github.com/leanprover-community/mathlib4/blob/b63493a4746b4651fceb3bcf3a4651cc36a4b8de/Mathlib/Data/Real/ConjExponents.lean#L59}{\extlink} from \mathlib.

We have now described most of the prerequisites for understanding the Lean statement of Theorem~\ref{metric-space-Carleson}, available as
\lean{metric_carleson}. The hypothesis \lean{hT} relies on the definition \lean{HasBoundedStrongType}\href{https://github.com/fpvandoorn/carleson/blob/AFM-frozen/Carleson/Defs.lean#L80-L86}{\extlink} to express the bound \eqref{nontanbound} for the
\lean{nontangentialOperator}.\href{https://github.com/fpvandoorn/carleson/blob/AFM-frozen/Carleson/Defs.lean#L211-L212}{\extlink}

\begin{lstlisting}[mathescape]
theorem metric_carleson {X : Type*} {a : ℕ} [MetricSpace X] {q q' : ℝ≥0}
  {F G : Set X} {K : X → X → ℂ} [KernelProofData a K] {f : X → ℂ}
  [IsCancellative X (defaultτ a)] (hq : q ∈ Ioc 1 2)
  (hqq' : q.HolderConjugate q') (mF : MeasurableSet F) (mG : MeasurableSet G)
  (mf : Measurable f) (nf : (‖f ·‖) ≤ F.indicator 1)
  (hT : HasBoundedStrongType (nontangentialOperator K · ·) 2 2 volume volume (C_Ts a)) :
  ∫⁻ x in G, carlesonOperator K f x ≤
    C1_0_2 a q * volume G ^ (q' : ℝ)⁻¹ * volume F ^ (q : ℝ)⁻¹
\end{lstlisting}

The \lean{carlesonOperator}\href{https://github.com/fpvandoorn/carleson/blob/AFM-frozen/Carleson/Defs.lean#L228-L229}{\extlink} is defined as a special case of the \lean{linearizedCarlesonOperator};\href{https://github.com/fpvandoorn/carleson/blob/AFM-frozen/Carleson/Defs.lean#L222-L224}{\extlink}
these definitions correspond to those introduced in equations \eqref{def-main-op} and \eqref{def-lin-main-op}, respectively. The expression \mmlean{‖carlesonOperatorIntegrand K (Q x) R₁ R₂ f x‖$_{\mathrm{e}}$} denotes the extended norm of the integral appearing in \eqref{def-T-r}, taking values in $\ennreal$ (see \S\ref{subsec:enorm} for more details).
In \lean{carlesonOperatorIntegrand}\href{https://github.com/fpvandoorn/carleson/blob/AFM-frozen/Carleson/Defs.lean#L216-L218}{\extlink} below, \lean{Annulus.oo x R₁ R₂}\href{https://github.com/fpvandoorn/carleson/blob/AFM-frozen/Carleson/Defs.lean#L55}{\extlink} denotes the set of points $y\in X$ such that $R_1 < \text{dist}(x, y) < R_2$.
\begin{lstlisting}[mathescape]
def carlesonOperatorIntegrand [FunctionDistances ℝ X] (K : X → X → ℂ)
    (θ : Θ X) (R₁ R₂ : ℝ) (f : X → ℂ) (x : X) : ℂ :=
  ∫ y in Annulus.oo x R₁ R₂, K x y * f y * exp (I * θ y)

def linearizedCarlesonOperator [FunctionDistances ℝ X] (Q : X → Θ X)
    (K : X → X → ℂ) (f : X → ℂ) (x : X) : ℝ≥0∞ :=
  $\sqcup$ (R₁ : ℝ) (R₂ : ℝ) (_ : 0 < R₁) (_ : R₁ < R₂),
    ‖carlesonOperatorIntegrand K (Q x) R₁ R₂ f x‖$_{\mathrm{e}}$

def carlesonOperator [FunctionDistances ℝ X] (K : X → X → ℂ) (f : X → ℂ)
    (x : X) : ℝ≥0∞ :=
  $\sqcup$ (θ : Θ X), linearizedCarlesonOperator (fun _ ↦ θ) K f x
\end{lstlisting}

Using these definitions, we can also state the theorem \lean{linearized_metric_carleson} (Theorem~\ref{linearized-metric-Carleson}). It requires an extra input in the form of a Borel function $Q\colon X \to \Theta$ with finite range, and its conclusion involves the
\lean{linearizedCarlesonOperator}\href{https://github.com/fpvandoorn/carleson/blob/AFM-frozen/Carleson/Defs.lean#L222-L224}{\extlink} $T_Q$. The hypothesis \lean{hT} encodes equation \eqref{linnontanbound}, written in terms of the
\lean{linearizedNontangentialOperator}.\href{https://github.com/fpvandoorn/carleson/blob/AFM-frozen/Carleson/Defs.lean#L204-L208}{\extlink}
\begin{lstlisting}[mathescape]
theorem linearized_metric_carleson {X : Type*} {a : ℕ} [MetricSpace X]
  {q q' : ℝ≥0} {F G : Set X} {K : X → X → ℂ} [KernelProofData a K]
  {Q : SimpleFunc X (Θ X)} {f : X → ℂ} [IsCancellative X (defaultτ a)]
  (hq : q ∈ Ioc 1 2) (hqq' : q.HolderConjugate q') (mF : MeasurableSet F)
  (mG : MeasurableSet G) (mf : Measurable f) (nf : (‖f ·‖) ≤ F.indicator 1)
  (hT : ∀ θ : Θ X, HasBoundedStrongType
    (linearizedNontangentialOperator Q θ K · ·) 2 2 volume volume (C_Ts a)) :
  ∫⁻ x in G, linearizedCarlesonOperator Q K f x ≤
    C1_0_2 a q * volume G ^ (q' : ℝ)⁻¹ * volume F ^ (q : ℝ)⁻¹ := by
\end{lstlisting}

\subsection{Verifying the definitions and statements}

The Lean kernel checks that all proofs are correct. It cannot ascertain that all definitions and theorem statements match their counterparts in the blueprint; this requires manual effort.
One strategy (going back to at least the Liquid Tensor Experiment \cite{LTEDefinitions})
is to prove basic results about these definitions.
In the best case, one can prove theorems that uniquely characterize a definition.
In our case, the classical Carleson theorem, \Cref{classical-carleson} is easy to state:
it uses besides \lean{partialFourierSum} only definitions in \mathlib.
This gives some evidence that the metric Carleson theorem is also correctly stated,
since it proved the classical theorem as a corollary.

Another aspect is making the relevant definitions easy to audit: to this end, we have collected all definitions and set-up necessary for stating the main theorems, as well as the main theorems' statements, into a single file. Auditing these 170 lines of Lean code (on top of \mathlib) ensures the theorems say what they claim --- over two orders of magnitude smaller than the entire project.
In addition, about a quarter of these lines are basic definitions that will be added to \mathlib soon, where they will receive careful review; this further reduces the need to trust code beyond \mathlib.

Lean ensures that the formalized proof is a valid proof of the given statement.
For maximum robustness (also suitable against adversarial code), we use the \emph{comparator} tool\footnote{\url{https://github.com/leanprover/comparator}} to verify that the final theorem we prove corresponds to the statement in the definitions file.
This verification is robust against malicious code or the abuse of Lean's metaprogramming facilities. It is the gold standard for verification of formalizations.

\section{Project organization}\label{sec:org}

Planning for the Carleson project began in Fall 2023 with
F.v.D.\ and the harmonic analysis group in Bonn,
consisting of L.B., Christoph Thiele and Rajula Srivastava.
At that time, L.B., Asgar Jamneshan, Rajula Srivastava and Christoph Thiele had written a 30-page proof of Theorem~\ref{metric-space-Carleson}, which turned into the preprint \cite{ThieleCarlesonPreprint}.
While the proof was precise and written with a lot of care,
this paper was intended for experts in harmonic analysis.
This meant that the proof omitted many arguments that are routine to experts,
but which would require quite some work to figure out for outsiders.

After a very brief attempt to formalize the paper directly,
we realized we needed a detailed blueprint that could be used as a foundation for the formalization.
The goal of the blueprint was that a non-expert in harmonic analysis
could read the statement and proof of any lemma in the blueprint,
along with all the definitions and lemmas that were explicitly referenced by that lemma,
and that alone would provide enough information for understanding the lemma and its proof, which is crucial in order to formalize it in Lean.
The precise way that this lemma fits into the rest of the proof was usually not specified,
and is not necessary for formalizing the proof of a single lemma.

The proof in the blueprint was divided into 179 lemmas.
To prepare for the collaborative formalization effort,
F.v.D.\ formalized the statements of most lemmas and the definitions referred to in them, which helped to ensure consistency across statements in the project. Since most contributors were only formalizing proofs, this also helped them to find the relevant definitions. To aid with this last point, especially for definitions only appearing inside proofs, we also manually maintained a table
that recorded the correspondence between notation and definitions in the blueprint on one side,
and Lean notation and declarations on the other.

During the public formalization effort, we used three platforms for communication and coordination,
an organizational structure that has also been used by other Lean projects,
such as the PNT$+$ project\footnote{\url{https://alexkontorovich.github.io/PrimeNumberTheoremAnd/}}
and the Equational Theories project \cite{tao2025equationaltheories}.
Firstly, we used GitHub to have a centralized version of the code,
using pull requests to ensure that at least one other contributor reviewed and approved additions or changes to the code.
Secondly, we had a webpage that contained an HTML version of the blueprint,
built using Patrick Massot's \texttt{leanblueprint} software.\footnote{\url{https://github.com/PatrickMassot/leanblueprint}}
This software generates a dependency graph, which gives a nice and quick overview of the progress of the project.
It also generates for each lemma in the blueprint convenient links to the corresponding Lean lemma(s).

Lastly, we used a public channel of Lean's Zulip chat to discuss the project and to post tasks needed for the formalization.
Most of these tasks consisted of formalizing a single lemma,
but some lemmas were split into multiple tasks, and some tasks consisted of something else,
like refactoring existing code in either the Carleson project or \mathlib.
We also used Zulip to keep track of the status of these tasks:
whether they were available, claimed or completed.
Zulip was also used in case someone had questions about a task,
or found a potential issue in the proof of a lemma.
In that case, another contributor or a member of Bonn's harmonic analysis group could help clarify the proof.
In a few cases, an inaccuracy was found that required changes to the statements of one or more lemmas in the blueprint,
as described in \S\ref{subsec:blueprint_mistakes}.

We divided the material to be formalized into two categories:
general material that would be useful in many other formalization projects,
and material specific to the current formalization, that would likely see little reuse.
We put the former material in a separate \texttt{ToMathlib} subfolder,
emphasizing the aim to upstream this material to \mathlib.
We also ensured as much as we could that these concepts and theorems would be developed in a greater generality than what was needed to formalize the current result, so that it would have an appropriate generality for \mathlib.

The material specific to this formalization was usually stated only at the required level of generality, and we relaxed some of the contribution requirements for this part of the project: for instance, we did not try to optimize all proofs or insist that they adhered to \mathlib's style conventions.
This helped make the Carleson project accessible to a broader range of potential contributors.

When the project was finished, it comprised 35.6k lines of code (excluding blank lines and comments), of which 9.6k were in the \texttt{ToMathlib} subfolder (using commit \texttt{\href{https://github.com/fpvandoorn/carleson/tree/32c2adc411c5c274ed8f19aa8041157e394fa38a}{32c2adc4}}).

\section{Working with a blueprint}\label{sec:blueprint}
\subsection{The blueprint writing process}\label{subsec:blueprint_writing}

As mentioned in \S\ref{sec:org}, the blueprint~\cite{CarlesonBlueprint} for the Carleson formalization project was written
by the harmonic analysis group in Bonn. When the formalization project was announced, it was not finished in its entirety: some parts of Section~7,
which deals with arguably the most complicated case of the proof of \Cref{discrete-Carleson},
were still being finalized, and
Section~3, which reduces \Cref{metric-space-Carleson} to \Cref{finitary-Carleson},
still had to be written in its entirety.
However, the rest of the blueprint was finalized
and contained more than enough material to start the formalization.

The main objective of the blueprint was to provide very detailed proofs for the lemmas and propositions
needed for the metric space Carleson theorem.
Another objective was to include a detailed proof
of the deduction of \Cref{classical-carleson} from \Cref{metric-space-Carleson}.
Finally, we wanted to streamline the argument to simplify the formalization.
One part of the formalization that we were worried about was the so-called ``invisible mathematics'' \cite{Bauer-2023} that one has to formalize.
These are arguments that are not written on paper,
but have to be provided by the formalizer to arrive at a formal proof.%
\footnote{There is also invisible mathematics that is performed by the proof assistant,
such as interpreting the meaning of notation or the scope of variables,
which we will not discuss here.}
Examples include proofs of integrability or measurability of side-goals, which are required
when manipulating expressions involving integrals.
To simplify these arguments, the blueprint was written to rely on finitary arguments.
The main part of the proof was done in a finite setting:
there are only finitely many cubes in the grid structure
and only finitely many tiles in the tile structure.
Moreover, we integrate over a subset of finite measure,
and the input function is a bounded function whose support has finite measure.
Only in the final steps of the proof,
multiple limiting arguments are used to derive \Cref{metric-space-Carleson}
from its finitary version, \Cref{finitary-Carleson}.
This approach also had downsides, as discussed in \S\ref{subsec:blueprint_lessons}.

\subsection{Changes and refinements}\label{subsec:blueprint_changes}

Throughout the formalization process, most of the blueprint only required minor changes, such as fixing typos or clarifying minor details. Let us highlight four aspects where formalization resulted in more substantial changes or refinements.

\subsubsection{$I \leq J$ vs.\ $I \subseteq J$}
The axioms for a grid structure are meant to convey that the scale of a cube roughly specifies the size of a cube.
Indeed, \eqref{eq-vol-sp-cube}, together with the fact that $D\gg 16$,
states that if we have cubes $I$ and $J$ with different scales, then they are subsets/supersets of balls with drastically different radii.
More precisely, if $I\subseteq J$ and $s(I)<s(J)$,
then $I$ is a subset of a ball with radius $4D^{s(I)}$
and $J$ is a superset of a ball with radius $\frac{1}{4}D^{s(J)}\gg 4D^{s(I)}$.
It is therefore tempting to conclude that $J$ must be a strict superset of $I$.
This might not be the case in an arbitrary metric space,
since there is no guarantee that there are points
in an annulus with specified inner and outer radii.

This also means that for two cubes $I$ and $J$, if $I\subseteq J$ as sets,
then we cannot conclude that $s(I)\le s(J)$.
The first version of the blueprint incorrectly described the data in a cube as just a subset of $X$,
and the scale function $s$ as a function on these subsets.
This implicitly assumed that $s(I)$ is only determined by $I$, viewed as a subset of $X$.
This was fixed by including the scale as part of the data of the cube.

In the blueprint we kept using the notation ``$I \subset J$'' for the inequality between cubes,
even though that actually included the assumption that $s(I) \le s(J)$.
This change had no further impact on downstream proofs, except for minor changes to lemma statements.

We thank Georges Gonthier for finding this inaccuracy in the blueprint.

\subsubsection{Just one top cube}

We made a second change to the definition of a grid structure, also suggested by Georges Gonthier.
The first definition did not include the axiom that there exists a top cube,
so it did not assume the existence of $I_0$ as in item~\ref{def:grid-structure}\ref{it:subsetmaxcube}.
In the construction of the grid structure on $X$
from Section~4 in the blueprint,
the top level always has a single cube,
so we were free to add this to the general definition of a grid structure.
This slightly simplified a few proofs in later sections,
and there was one computation that already implicitly assumed
that there was only one top cube.

\subsubsection{The Hardy--Littlewood maximal function} An important ingredient in the proof is a Vitali covering argument, using the Hardy--Littlewood maximal function.
Recall that, given $f\colon\real^d\to\complex$, the (centered) \emph{Hardy–Littlewood maximal function}\href{https://github.com/fpvandoorn/carleson/blob/AFM-frozen/Carleson/ToMathlib/HardyLittlewood.lean#L31-L32}{\extlink} $M f$ of $f$ is defined as
\begin{equation*}
Mf\colon\real^d \to\ennreal, \quad x\mapsto \sup_{r>0}\frac{1}{|B(x,r)|}\int_{B(x,r)}|f(y)|\,dy.
\end{equation*}
If $f$ is integrable, then $M f$ is (weak) $L^1$, and hence finite almost everywhere.\href{https://github.com/fpvandoorn/carleson/blob/AFM-frozen/Carleson/ToMathlib/HardyLittlewood.lean#L365-L369}{\extlink} Initially, the blueprint only used a variant taking a supremum over finitely many balls: this ensured the resulting function was \emph{always} finite, hence real-valued. Applying it to a countable collection of balls requires taking a limit and using the monotone convergence theorem.
This is somewhat inelegant, since it uses a non-standard definition and requires an extra limiting argument.

For the formal proof, we used the standard definition\href{https://github.com/fpvandoorn/carleson/blob/AFM-frozen/Carleson/ToMathlib/HardyLittlewood.lean#L31-L32}{\extlink}
using arbitrary suprema: this means the maximal function is $\ennreal$-valued.
Being able to speak of $L^p$-functions valued in $\ennreal$ prompted the introduction of \emph{extended norms}, see \S\ref{subsec:enorm}.
We also generalized the allowed input functions to any function from a pseudo-metric space to any space with an extended norm.

\subsubsection{The real interpolation theorem}\label{subsubsec:real_interpolation}
The Marcinkiewicz real interpolation theorem is a standard result in functional analysis.
It is also a key prerequisite for the proof of Carleson's theorem: for example, it is applied to the Hardy–Littlewood maximal function to deduce it is of strong type $(p, p)$, for $p\in(1,\infty)$.
Our formalized version aims to reduce superfluous assumptions commonly found in the mathematical literature. For instance, it applies to any sub-additive operator (not merely sub-linear ones) --- that is, there is no need for a scalar multiplication on the target. Similarly, any positive exponent $p$ is allowed (as opposed to demanding $p\geq 1$).

In light of generalizing the Hardy--Littlewood maximal function, this proof had to be generalized again:
applying it to the maximal functions requires a version which applies to functions with co-domain $\ennreal$ (or, more generally, any monoid with a continuous extended norm function).
In return for a significant amount of refactoring work,\footnote{The initial formalization was about 5000 lines long; the generalization touched about 2000 of them.}
we have formalized a more general argument that is also conceptually clearer: for instance, it demonstrates that the subtractive structure on the co-domain is not used in a meaningful way. We ended up with the following version.
\begin{theorem}[real interpolation theorem\href{https://github.com/fpvandoorn/carleson/blob/AFM-frozen/Carleson/ToMathlib/RealInterpolation/Main.lean\#L1439-L1448}{\extlink}]
    Let $(X,\mu)$ and $(Y,\nu)$ be measure spaces, let $E$ be a topological space with the structure of an additive monoid endowed with a compatible continuous extended seminorm, let $F$ be a topological space with a continuous extended norm function.
    Let $p_0, p_1, q_0, q_1, p, q \in \ennreal$ such that $p_0 \in (0,q_0]$ and $p_1 \in (0,q_1]$ and $q_0 \ne q_1$ and there exists $t \in (0,1)$ such that $\frac{1}{p} = \frac{1-t}{p_0} + \frac{t}{p_1}$ and $\frac{1}{q} = \frac{1-t}{q_0} + \frac{t}{q_1}$.
    Let $T\colon (X \to E) \to (Y \to F)$ be an operator that is subadditive in the following sense: there exists a constant $A \ge 1$ such that for $f, g \in L^{p_0}(X; E) \cup L^{p_1}(X; E)$, we have $\|T(f+g)(y)\| \le A \cdot \left(\|Tf(y)\| + \|Tg(y)\|\right)$ for $\nu$-a.e.\ $y\in Y$. Furthermore, we assume that for any $f\in L^p(X;E)$ the function $T f$ is $\nu$-a.e.\ strongly measurable.

Then, if $T$ has both weak type $(p_0, q_0)$ and weak type $(p_1, q_1)$, it has strong type $(p,q)$.
\end{theorem}

More details will be provided in an upcoming article \cite{PortegiesRothgangRealInterpolation}.

\subsubsection{Constant tweaking}\label{subsubsec:const-tweak}
As explained below in
\S\ref{subsec:constants}, the constants appearing in the blueprint
are always explicit and most often in the form $2^{n a^3}$ where $a$ is a
given parameter of the system and $n$ is a concrete natural number, ranging
from $1$ or $2$ in the first lemmas to $443$ in the final result.

When a result depends on several previous lemmas, its constant can typically
be expressed in terms of the previous lemmas, giving chains of conditions for
the overall validity of the constants. When correcting a minor mistake in the
proof of a lemma requires a slight increase in its constant, all statements
that depend transitively on it must also be corrected. The standard practice
in paper mathematics to avoid this kind of intricate dependency is not to
care about the constants, and formulate the lemmas as ``there exists a
universal constant $C$ such that ...''. This is not a viable approach in
formalized mathematics, as explained in \S\ref{subsec:constants}.

To mitigate this issue and avoid endless corrections, during the development
process the constant in each lemma did not have a numerical value, but
instead it was expressed as a function of the constants in the other lemmas.
Therefore, changing the constant in Lemma~A would automatically change the
constant of Lemma~B that depends on Lemma~A, not altering the validity of its
proof. At the very end of the formalization process, when all the proofs
worked, we could finally register the numerical values of the constants with
a correctness guarantee, starting with the first lemmas and propagating the
values iteratively. The blueprint was then fixed using these computer-checked
constants --- and indeed many constants written in earlier versions of the
blueprint were changed, although this had absolutely no consequence on the validity of
the proof strategy.

\subsection{Dealing with inaccuracies}\label{subsec:blueprint_mistakes}

Our formalization of the metric space Carleson theorem certified that the main result in \cite{ThieleCarlesonPreprint} was correct, and that its detailed proof provided in \cite{CarlesonBlueprint} did not contain any serious mathematical mistakes. However, the formalization surfaced various small inaccuracies in the blueprint.  This is not surprising, since experience shows that a human-written mathematical document of this scale usually contains at least a few minor technical flaws.
Most of these flaws were not present in the original paper; instead they were introduced during the blueprint writing process, which required providing a much higher level of detail than in a standard harmonic analysis paper.

Many proofs in the blueprint relied on non-standard finitary arguments, which sometimes led to proofs that were not correct in the edge cases.
In some places in the blueprint, proofs were slightly inaccurate, and at other places lemmas required stronger assumptions than originally stated. However, when assumptions had to be strengthened, the
stronger assumptions were always satisfied at the places where the lemmas
were used. Therefore, all these issues were purely local and did not impact the global proof.

In most cases, the mistakes detected in the blueprint were minor and could be resolved by the formalizers themselves, for instance by adding a missing hypothesis or by modifying a certain constant. However, fixing certain proofs required the assistance of harmonic analysis experts with a deep understanding of the overall proof strategy. For instance, the proof of an earlier version of \cite[Lemma~6.2.3]{CarlesonBlueprint} relied on the fact that a certain pair of dyadic cubes were not disjoint, which was not necessarily the case under the lemma assumptions. The non-obvious solution required adding an extra hypothesis that a certain pair of open balls were not disjoint and increasing a constant appearing in the lemma statement, after which the previous proof could be adapted to conclude the result. Luckily this error did not propagate, since the only application of the lemma was still possible under the new hypotheses.

A more involved fix was required for \cite[Lemma~6.3.4]{CarlesonBlueprint}, whose original proof presented several issues. The first needed change was simple, and consisted of further restricting the allowed values for the parameter $\vartheta$, so that a certain covering property could be applied. However, the proof was also incorrect in the case where the value $s(\mathfrak{p})$ associated to a certain tile $\mathfrak{p}$ took the minimum possible value, and it incorrectly assumed that a certain new tile constructed during the proof would always be distinct from $\mathfrak{p}$. To repair these issues, the statements of both Lemma~6.3.3 and Lemma~6.3.4 had to be modified, and significant changes had to be applied to the proof of Lemma~6.3.4.

Another problem was found in Lemma~11.1.6 of the blueprint, which is a version of the boundedness of the
Hilbert transform on $L^2(\real)$. It is proved by localization and
periodization, ultimately using that $\|S_N f\|_{L^2[0,2\pi]} \leq
\|f\|_{L^2[0,2\pi]}$ for periodic functions. One step in the proof
involves a convolution $\int f(y) g(x-y) \ud y$ for some kernel $g$. In
the initial argument, the position of $x-y$ was not carefully
controlled, claiming to use a universal bound on $g$ while this bound
was only valid in some limited interval. The issue was fixed by
truncating $f$ and only looking at $x$ in some interval to make sure
that, when $f(y)\ne 0$, then $x-y$ belonged to the place where $g$ is
well controlled. These additional assumptions were satisfied in the
rest of the proof.

Finally, Lemma~8.0.1 in the blueprint is a regularization lemma, saying that
some Hölder functions can be approximated by Lipschitz functions,
constructed through a smoothing integral formula. The original proof
assumed implicitly that a continuous function which is
$\tau$--Hölder--continuous in a closed ball $B(z,R)$ and supported
there is globally $\tau$--Hölder--continuous with the same control.
This is true in well-behaved spaces, for instance in $\real^n$, but not
in the generality of the paper. For instance, if the annulus $\{y \mid
d(z,y) \in [R-\epsilon, R+\epsilon]\}$ is empty, then the
characteristic function of $B(z, R)$ is Hölder continuous on the ball
with a uniform bound and supported there, but its extension to the
whole space has a norm which may grow like $\epsilon^{-\tau}$. The fix
was to modify the assumptions of the lemma, assuming from the start
that the function is $\tau$--Hölder--continuous in the larger ball
$B(z, 2R)$. This more restrictive assumption is satisfied in all the
places where the lemma is used, although this forced us to alter
some constants in other lemma statements.

\subsection{Lessons learned}\label{subsec:blueprint_lessons}
In hindsight, some slightly different design choices during
the writing of the blueprint could have further streamlined the formalization process.

One of these was the use of finitary arguments, specifically the non-standard
definition of the Hardy--Littlewood maximal function over finite collections of balls,
and the decision to only allow finitely many scales in the constructions of grid
and tile structures.

While these choices were made partly with the intent to simplify, they ultimately
led to additional technical complications in the blueprint, causing several inaccuracies
and some friction during the formalization.

In the case of the maximal function, the finitary formulation required explicitly specifying
various finite collections of balls, some of which were initially slightly off. As mentioned
above, the finitary variant also led to an unnecessary limiting argument to generalize
to the standard version.

The blueprint required that all dyadic cubes in a grid structure (\Cref{def:grid-structure}) have a scale in a fixed interval $[-S, S]$.
This requirement is not present in the original preprint~\cite{ThieleCarlesonPreprint}.
While this finiteness assumption simplified parts of the argument, some inaccuracies were introduced in the blueprint related to arguments at the bottom scale.

Some general results, such as the Marcinkiewicz real interpolation theorem, were proven in the blueprint only in the generality required to apply them in the proof of the main result. However, we ended up formalizing a more general version so that these theorems could be reused in other formalization projects. This would have been easier if their statements and proofs had already been generalized during the writing of the blueprint.

A final point is that the lack of definition environments in the blueprint also
meant that definitions were absent from the dependency graph. This was especially problematic for those contributors formalizing statements, since sometimes the definitions appearing in them had been declared a few sections before and there was no convenient way to find them. This issue was mitigated by maintaining a table for the definitions, as explained in \S\ref{sec:org}. This table was sparse at the beginning of the formalization, but more definitions were added to it as the project advanced.

\section{Design decisions}\label{sec:design}
\subsection{Treatment of constants}\label{subsec:constants}
In analysis papers, one often needs to compare the growth of two functions $f, g\colon X \to Y$, where $X, Y$ are real or complex normed vector spaces. This is often expressed by statements of the form
\begin{align}\label{eq:const}
	\exists C > 0, \forall x \in X, \| f (x) \| \le C\|g(x)\|.
\end{align}
In many cases, one only cares about the existence of the positive constant $C$, but not about its concrete value.
The constant is often left implicit, by writing \eqref{eq:const} as
\[ \forall x \in X, \| f (x) \| \lesssim \|g(x)\|.\]

This was also the case in the first version of \cite{ThieleCarlesonPreprint}. However, while Lean allows one to work with existential statements such as \eqref{eq:const}, this would have made the formalization harder. From the mathematical point of view, it would require the formalizers to figure out how large these constants needed to be and how the constants required to make each lemma hold were related to each other. From the technical side, it would require working with witnesses for these existential statements, which can introduce some overhead.
Therefore, while they were writing the blueprint, the authors of \cite{ThieleCarlesonPreprint} decided to make all of the constants explicit. At the beginning of the blueprint, they fixed a real number $1 < q \le 2$ and a natural number $a \ge 4$, which is used to define the constants $D = 2^{100a^2}, \kappa := 2 ^{-10a}$ and $Z := 2^{12a}$. The bounding constants that appear in the blueprint theorem statements are expressed as functions of these fundamental constants.

In the formalization, we introduced a constant \lean{c} with a value of 100, coming from the exponent of $D$.
Every other constant is defined in terms of $a$, $q$ and \lean{c}, and we use the notation convention \lean{Cx_y_z} to denote the constant appearing in result \lean{x.y.z} of the blueprint.

Once the main result had been formalized, we were able to check that we could replace the constant $D := 2^{100a^2}$ by $2^{7a^2}$ and the results still hold without needing to make any adjustments to the proofs. This improves the constant $2^{443a^3}$ appearing in \Cref{metric-space-Carleson} to $2^{46a^3}$.
This constant could be improved further by adding more fine-grained constants,
instead of having a single constant $a$ that plays multiple roles.

\subsection{Standing assumptions and the \texttt{ProofData} pattern}\label{subsec:proof_data}

Mathematical papers often contain standing assumptions, that are in place for
all the statements in the paper or in a specific section. This avoids
tedious repetition in the statements, sometimes at a slight expense in
readability as the reader has to locate the places where such standing
assumptions are made. Such a practice is especially important when the
standing assumptions are lengthy and technical.

The standing assumptions for the statement of the generalized
Carleson theorem are quite formidable,
see~\Cref{metric-space-Carleson}.
They involve a metric space $X$, a measure $\mu$ on $X$ which is
doubling for a constant $2^a$, a family $\Theta$ of functions from $X$ to
$\real$, a distance on $\Theta$ satisfying five properties as well as a
cancellativity condition, and a kernel $K\colon X \times X \to \real$ with
Calderón--Zygmund like properties. The proofs involve even more additional
data, notably two bounded measurable sets $F$ and $G$, an integrability.
In the proof we introduce additional data,
notably two bounded measurable sets $F$ and $G$, an integrability
exponent $q\in (1, 2]$ and two measurable functions $\sigma_1$ and
$\sigma_2$ with finite range.

One possibility in the formalization would be to include
all of these as arguments to each lemma that uses them.
This would come at a high readability cost. Moreover, each time
one applies such a lemma one would need to provide the individual assumptions
as a long list of parameters to the lemma.

It is better to try to also use the practice of
standing assumptions in the formalization, in some form. One could define a
structure containing all these standing assumptions, and have it as a
parameter in every lemma and every definition. The readability cost would be
minor, but still present. In the Carleson project, we have used another
strategy, making all the standing assumptions completely implicit --- in
particular, we do not need to pass them between lemmas. The implementation
takes advantage of a feature of the language, type classes, initially designed
for a completely different purpose.

The way to declare that a space is a metric space is the following line:
\begin{lstlisting}[mathescape]
variable {X : Type*} [MetricSpace X]
\end{lstlisting}

The declaration between square brackets is the \emph{type class}. When one writes the
distance \texttt{dist x y} between two points $x$ and $y$ of $X$, the system
tries to locate such a metric space type-class instance on $X$, and uses the distance
provided by this instance. The same system is used to provide algebraic
operations on types, and properties of these.

Type classes can be used to capture standing assumptions as follows.
We declared a new type class, of the form \lean{[KernelProofData X]},\href{https://github.com/fpvandoorn/carleson/blob/AFM-frozen/Carleson/Defs.lean#L261-L266}{\extlink}
containing a measure on $X$,
a doubling constant, a family $\Theta$ of functions from $X$ to
$\real$, and so on, i.e., all our standing assumptions. Then, whenever one
tries to use an object from the standing assumptions, the system looks for a
\lean{KernelProofData} instance and fetches the object from the type-class
instance, just like it fetches the distance in a metric space.

Another advantage of this approach is that type classes can extend each other
(just like a normed space contains more information than a metric space).
This means we can craft one type class containing the standing assumptions for
the statements, and another extended type class \lean{ProofData}\href{https://github.com/fpvandoorn/carleson/blob/AFM-frozen/Carleson/ProofData.lean#L20-L36}{\extlink} with more
objects for the proofs. When a \lean{ProofData} instance is in scope but
Lean needs a \lean{KernelProofData} instance, type-class inference
automatically fetches the latter from the former, just as it transparently
interprets a normed space as a metric space.

The specific implementation of this idea in the Carleson project is slightly
more involved, as it is written in the form \lean{[ProofData a q K σ₁ σ₂ F G]}. However, all the additional arguments $a$, $q$, $K$, $\sigma_1$,
$\sigma_2$, $F$ and $G$ are inferred automatically from the selected type-class instance, and the user only has to specify the type $X$. In
technical terms, all arguments other than $X$ are \emph{output parameters}.
The user does not have to specify the output parameters explicitly,
as expected of a system designed to
handle standing assumptions.

\subsection{Working with real numbers}\label{subsec:real}

The project differentiated between three different types related to the real numbers, $\real$,\href{https://github.com/leanprover-community/mathlib4/blob/b63493a4746b4651fceb3bcf3a4651cc36a4b8de/Mathlib/Data/Real/Basic.lean#L36-L38}{\extlink} $\nnreal$ \href{https://github.com/leanprover-community/mathlib4/blob/b63493a4746b4651fceb3bcf3a4651cc36a4b8de/Mathlib/Data/NNReal/Defs.lean#L58}{\extlink}
and $\ennreal=[0,\infty]$,\href{https://github.com/leanprover-community/mathlib4/blob/b63493a4746b4651fceb3bcf3a4651cc36a4b8de/Mathlib/Data/ENNReal/Basic.lean#L101}{\extlink}
which created difficulties in the formalization
and tension in design decisions.
Indeed, which of the three types should be used for constants,
suprema, integrands, and computations?

Since every term has a unique type, coercions can be registered,
and afterwards are automatically inserted to go from one type to another. In particular, there are coercions from $\nnreal$ to both $\real$ and $\ennreal$.
In proofs, working with coercions is somewhat annoying, but the \texttt{norm\_cast} tactic is useful to normalize such casts.

However, one also needs conversions in the other directions, for instance to convert real numbers to non-negative real numbers. These conversions commute with other operations like addition and multiplication, but only under certain conditions. For example, since \texttt{Real.toNNReal}\href{https://github.com/leanprover-community/mathlib4/blob/b63493a4746b4651fceb3bcf3a4651cc36a4b8de/Mathlib/Data/NNReal/Defs.lean#L154-L156}{\extlink} maps negative numbers to zero, we get the following inequality.
\begin{lstlisting}
Real.toNNReal (-2) * Real.toNNReal (-3) ≠ Real.toNNReal ((-2) * (-3))
\end{lstlisting}
Even with the necessary side conditions, the \texttt{norm\_cast} tactic does not get very far
in these situations.
This can make conversions between number types quite painful,
and adds some pressure to limit conversions and put as many numbers as possible in the same number type.

Sometimes one is forced to add a conversion, because certain operations are not defined on all types.
Consider, for instance, performing a computation on exponents in $L^p$-function spaces.
Around the real interpolation theorem, one needs to compute $p$ in terms of other exponents that
naturally live in $\ennreal$.
Yet there is currently no definition for $a^p$ if $a$ and $p$ are in $\ennreal$,
and therefore we used conversions from $\ennreal$ to $\real$.

There are also advantages of differentiating between these number types.
An advantage of putting for instance a constant in \texttt{NNReal} over \texttt{ENNReal} or
\texttt{Real}, is that it encodes that the constant cannot be infinity and that it needs to be non-negative.

Another advantage of differentiating comes from the different levels of tactic support.
There is better tactic support for $\real$ than for $\nnreal$ and $\ennreal$:
tactics such as \texttt{linarith} and \texttt{field\_simp} can simplify computations in $\real$.
For terms in $\nnreal$, support in \texttt{linarith} was only added after this project was completed.\footnote{\url{https://github.com/leanprover-community/mathlib4/pull/35155}}
Neither \texttt{linarith} nor \texttt{field\_simp} support terms in $\ennreal$.
The tactic \texttt{ring} \emph{does} work for these types, solving goals that only contain addition, multiplication and numerals.
In both $\ennreal$ and $\nnreal$, it will generally fail when a goal needs properties of subtraction: e.g.\ it fails to prove that $x - x = 0$.
In $\ennreal$, it also fails to prove equalities related to division or multiplicative inverses. If $x\colon \nnreal$, the \texttt{ring} tactic can prove that $(x^{-1})^2 = (x ^ 2)^{-1}$ and $3 / x = 3 x^{-1}$, but this fails if $x\colon \ennreal$.
Similarly \texttt{norm\_num} also works for numbers in $\ennreal$, but again not as well as for real numbers. Both \texttt{ring} and \texttt{norm\_num} fail to solve $(3/5) * (5/3) = 1$ when these numbers are viewed as numbers in $\ennreal$.

On the nose, many results about real numbers transfer to $\ennreal$ provided all arguments are finite; in many cases, these finiteness hypotheses are actually necessary. Manipulating expressions in $\ennreal$ thus often requires proving such finiteness side goals. The \lean{finiteness}\href{https://github.com/leanprover-community/mathlib4/blob/b63493a4746b4651fceb3bcf3a4651cc36a4b8de/Mathlib/Tactic/Finiteness.lean#L73}{\extlink} tactic (initially developed for the formalization of the PFR project\footnote{\url{https://teorth.github.io/pfr/}}) solves many such proof obligations automatically. We made this tactic more powerful, so that it supports all relevant cases from the Carleson project.

A desirable future step is adding improved interaction with the \lean{positivity}\href{https://github.com/leanprover-community/mathlib4/blob/b63493a4746b4651fceb3bcf3a4651cc36a4b8de/Mathlib/Tactic/Positivity/Core.lean#L593}{\extlink} tactic (which solves goals of the form \texttt{$0 \leq x$} or \texttt{$0 < x$}). Proving a goal \texttt{$(x + y + 1)^{-1}<\infty$} (for $x,y\in\ennreal$) requires proving that $x + y + 1$ is non-zero; this is already handled by calling \lean{positivity}. Conversely, \lean{positivity} alone cannot prove goals \texttt{$x^{-1}\neq 0$} for $x\colon\ennreal$: since this is only true when $x$ is finite, this would require that \lean{positivity} calls \lean{finiteness}. This is possible, but is not yet used for division or multiplicative inverses of numbers in $\ennreal$.

A more powerful tactic that would have been useful in the Carleson project is a tactic that can reduce goals in $\ennreal$ to $\real$: this could be implemented by case-splitting on each relevant variable being finite or infinite. The finite case allows reduction to $\real$ arithmetic; the infinite case will usually be a very simple calculation (as, for example, $\infty$ added to any number in $\ennreal$ is still $\infty$). Such goals can hopefully be proven automatically, or reduced to much simpler tasks for the user.
This was implemented in Mathlib in August 2026,\href{https://github.com/leanprover-community/mathlib4/pull/43254}{\extlink} after the completion of the Carleson project.

The choice of real number type also came up when integrating over a subset of real numbers. Before the Carleson project, there was no definition of integration over \lean{NNReal} or \lean{ENNReal}. As integrating such functions was useful, we added \lean{MeasureSpace}\href{https://github.com/leanprover-community/mathlib4/blob/b63493a4746b4651fceb3bcf3a4651cc36a4b8de/Mathlib/MeasureTheory/Measure/MeasureSpaceDef.lean#L355-L356}{\extlink} instances on these types and proved basic properties and relations. Now, the following three expressions to integrate a function $f : \ennreal\to \ennreal$ are defined and proven to be equal.
\begin{lstlisting}[mathescape]
∫⁻ (x : ℝ) in Ioi 0, f (ENNReal.ofReal x)
∫⁻ (x : ℝ≥0), f (ENNReal.ofNNReal x)
∫⁻ (x : ℝ≥0∞), f x
\end{lstlisting}
This avoids another source of casting between real number types.

As a final example, consider the supremum of a set.
In the Carleson project, we mostly used suprema where we knew the supremum was finite,
so we could consider a setup where we take suprema over subsets of real numbers
and the suprema are real numbers as well.
This was the choice at the outset of the project.
However, in this setup, we need to constantly prove that the set is bounded
or carry around the appropriate assumptions.
This is cumbersome, and therefore it turned out to be nicer to
work with suprema in $\ennreal$.

Gradually, we came to the conclusion that, with some exceptions, it was best to use $\ennreal$ as much as possible, which is also exemplified by the next section.

\subsection{Use of \texttt{ENorm}}\label{subsec:enorm}

This project motivated the creation of a new formalization abstraction: \emph{extended norms} (\lean{enorm}s for short) generalize many results from normed spaces to $\ennreal$.
Their introduction was motivated by defining the Hardy--Littlewood maximal function with codomain $\ennreal$ (see \S\ref{subsec:blueprint_changes}).
Previously, a statement ``$f\in L^p$'' was only defined for functions $f\colon X\to E$ from a measure space $X$ to a normed space $E$. However, speaking about the maximal function being in weak $L^p$ requires a notion of $L^p$ functions with target $\ennreal$. Fortunately, many basic facts about $L^p$ functions do not require the target to be a normed space and have reasonable analogues for $\ennreal$: for instance, the $L^p$-norm of $f\colon X\to\ennreal$ for $0<p<\infty$ can be defined as $(\int_X f(x)^p \ud x)^{1/p}$, where the integrand is the function $X\to \ennreal$ with $x\mapsto f(x)^p$.

Motivated by this observation, we introduced a new notation class \lean{ENorm},\href{https://github.com/leanprover-community/mathlib4/blob/b63493a4746b4651fceb3bcf3a4651cc36a4b8de/Mathlib/Analysis/Normed/Group/Defs.lean#L70-L74}{\extlink} describing a type endowed with a function $\texttt{enorm}\colon X\to\ennreal$.
This captures both normed spaces and $\ennreal$. The enorm of a normed space is the ambient norm, considered as a function into $\ennreal$; the enorm on $\ennreal$ is the identity function.
Some definitions only require the existence of an enorm; many results require suitable compatibility conditions.
Most lemmas necessary in this project only require an \lean{ESeminormedAddMonoid},\href{https://github.com/leanprover-community/mathlib4/blob/b63493a4746b4651fceb3bcf3a4651cc36a4b8de/Mathlib/Analysis/Normed/Group/Defs.lean#L111-L117}{\extlink} i.e.\ a monoid endowed with a continuous enorm that is positive semi-definite and satisfies the triangle inequality. Many lemmas, in fact, only require slightly weaker assumptions.

Pursuing this approach required generalizing many basic measure-theoretic definitions from \mathlib to \lean{enorm}s, as well as all basic results that are required in this project. This was a significant effort, but the verified nature of formalization helped a lot by allowing us to do so incrementally: in a first step, the \lean{ENorm} class was added, together with the enorm instance on normed spaces. As the second step, all basic definitions were changed to allow an \lean{enorm} as codomain, whenever this did not require further changes. The third step modified many lemma statements to use \lean{enorm}s instead of norms, without yet generalizing any type classes. The fourth step consisted of generalizing lemmas, by changing typeclasses such as \lean{NormedAddCommGroup} to \lean{ENormedAddCommMonoid}. Often, this required small changes to the lemma statements and proofs, such as arguments changing from a real to an extended real number, or adding a hypothesis about some argument in $\ennreal$ being finite. Lean helped with the refactor by ensuring correctness of all the changes.

To illustrate the magnitude of this refactoring, let us give some statistics: step two (changing the basic definitions) changed about 300 lines of code; step three (changing lemma statements to \lean{enorm}s) about 1000, and step 4 (lemma generalizations) touched about 2500 lines of code.

The last step is not exhaustive: while substantial parts of \mathlib have been generalized, some areas would require further refactor (and will follow as needed). An example of a refactor\href{https://github.com/leanprover-community/mathlib4/pull/42662}{\extlink} that happened after the Carleson project was completed was to generalize the total variation of a function enough to apply to functions with $\ennreal$ as codomain. This required generalizing the concept of metric space, since in \mathlib the type $\ennreal$ is not an (extended) metric space. Extended metric spaces require that the topology is generated from the distance, which is not the case in $\ennreal$.
Prior to this change, a weaker notion of extended metric spaces was introduced to \mathlib.\href{https://github.com/leanprover-community/mathlib4/pull/38105}{\extlink} Both of these refactors were performed by Felix Pernegger.

The process of generalization revealed that this abstraction is also useful on its own: previously, many proofs in \mathlib involved converting between norms taking values in $\real$ or $\nnreal$; using \lean{enorm}s avoids these conversions and is a natural way to phrase these statements.

\subsection{Working with \texorpdfstring{$L^p$}{L\^{}p} functions}\label{subsec:Lp}

Working with $L^p$ functions presents a further design choice: one option is to work with elements of $L^p$ space
\lean{MeasureTheory.Lp},\href{https://github.com/leanprover-community/mathlib4/blob/b63493a4746b4651fceb3bcf3a4651cc36a4b8de/Mathlib/MeasureTheory/Function/LpSpace/Basic.lean#L88-L96}{\extlink}
which are equivalence classes of functions up to almost everywhere equality;
a second option is to work with explicit functions and the condition that they are (a.e.\ strongly) measurable and that the $p$-th power of their norm is integrable, which is defined by the predicate \lean{MemLp}.\href{https://github.com/leanprover-community/mathlib4/blob/b63493a4746b4651fceb3bcf3a4651cc36a4b8de/Mathlib/MeasureTheory/Function/LpSeminorm/Defs.lean#L116-L119}{\extlink}

The mathematical difference is very small: virtually all operations in the project respect a.e.\ equality.
Formalization ergonomics, however, are greatly different: working with equivalence classes is much more cumbersome, as identities on the level of representatives only hold almost everywhere.
For example, $(f+g)(x) = f(x) + g(x)$ holds on the nose for actual functions.
However, if we have elements $\tilde{f},\tilde{g}\in L^p$,
then we can still write down $\tilde{f}(x)$ by picking an arbitrary representative of this class, but then $(\tilde{f}+\tilde{g})(x)=\tilde{f}(x)+\tilde{g}(x)$ only holds for almost every $x$.
This can make it significantly more cumbersome to manipulate expressions involving elements of the $L^p$ space, so it is easier to directly work with functions with the predicate \lean{MemLp}.
In simple examples, the \lean{filter_upwards}\href{https://github.com/leanprover-community/mathlib4/blob/b63493a4746b4651fceb3bcf3a4651cc36a4b8de/Mathlib/Order/Filter/Defs.lean#L444}{\extlink} tactic
allows us to conveniently use a.e.-equalities in proofs.
However, in more complicated examples, equality proofs using elements of $L^p$ get cumbersome.
The equation below is true by reflexivity if \lean{f} and \lean{g} are functions \lean{ℝ → ℝ},
but requires a tedious proof if they are elements of the $L^p$ space.
\begin{lstlisting}
∫ x, ∫ y, (f + g) (x * y) = ∫ x, ∫ y, f (x * y) + g (x * y)
\end{lstlisting}
For this reason, lemmas in \mathlib are mostly stated for \lean{MemLp} instead of \lean{Lp}.

This is analogous to the trade-off between bundled and unbundled objects.
In \mathlib various functions properties are defined in an unbundled way, for example a continuous function $f\colon X \to Y$ is usually formalized as \lstinline|{f : X → Y} (hf : Continuous f)| instead of using the bundled \lean{f : ContinuousMap X Y}.\href{https://github.com/leanprover-community/mathlib4/blob/b63493a4746b4651fceb3bcf3a4651cc36a4b8de/Mathlib/Topology/ContinuousMap/Defs.lean#L33-L37}{\extlink}

This relates well with the use of \lean{enorm}s (\S\ref{subsec:enorm}): \mmlean{Lp E p $\mu$} as a normed space requires functions with codomain in a Banach space; this excludes e.g.\ $\ennreal$.
Phrasing lemmas using \lean{MemLp} allows us to generalize normed spaces to spaces with an \lean{enorm}.

\subsection{Test function classes}\label{subsec:bdd_comp_supp}
A common task during the formalization was to prove the various integrability and measurability side goals
arising when estimating integrals.
As is typical in harmonic analysis, most technical estimates are {\it a priori} estimates involving
input functions in suitable test function classes.

In many parts of this project, bounded measurable compactly supported functions provided a convenient
test function class. We also used the slightly larger class
of bounded measurable functions supported on a set of finite measure in places where this weaker hypothesis
matches the blueprint more closely.
Since most constructions appearing in the proof leave these function classes invariant,
integrability and measurability considerations are usually mathematically trivial.

It is desirable for this simplicity to be reflected in the formalization.
This motivated packaging the recurring side conditions into test function classes.
The auxiliary structures \lean{BoundedCompactSupport}\href{https://github.com/fpvandoorn/carleson/blob/AFM-frozen/Carleson/ToMathlib/BoundedCompactSupport.lean#L84-L87}{\extlink}
and \lean{BoundedFiniteSupport}\href{https://github.com/fpvandoorn/carleson/blob/AFM-frozen/Carleson/Defs.lean#L74-L76}{\extlink}
contain bounded measurable functions with compact support and with support of finite measure, respectively.

Once these hypotheses are bundled, later proofs can invoke short API lemmas instead of repeatedly
reproving the same measure-theoretic facts. For example, from a hypothesis
\lean{hf : BoundedCompactSupport f} one immediately obtains that \lean{f} is a member of
any desired $L^p$ space,
or that \lean{s.indicator f} again has bounded compact support when \lean{s} is measurable.
Lemmas such as \lean{hf.carlesonSum}\href{https://github.com/fpvandoorn/carleson/blob/AFM-frozen/Carleson/Operators.lean#L164-L165}{\extlink} show that the relevant
operators preserve the test function class, so once compact support has been established at the input, it remains available throughout the argument.

Choosing the right definition of these structures can be subtle. With the
existing measure theory developments in \mathlib in mind, one might be tempted
to use the notion of {\em essential boundedness} (i.e.\ boundedness
potentially up to a set of measure zero) rather than {\em everywhere
boundedness} (and similarly, one might want to use {\em essential
compactness} rather than genuine compactness).

In this context, it seems preferable to impose the stronger requirements. In analysis, test functions
are typically chosen to have strong properties, making it easier to prove {\em a priori} estimates,
which can then be extended via density arguments.
Bounded functions of compact support are a natural function class that stays invariant with respect
to many common operations in analysis such as convolutions, fiberwise integration and projection.
Weakening assumptions to almost everywhere variants destroys this invariance for operations
related to measure zero sets.
For example, if $g$ is bounded with compact support, and $f$ is a closed embedding,
is $g\circ f$ a bounded function with compact support?
This is only true if $g$ has compact support and is bounded everywhere;
otherwise $f$ may send a set of positive measure to a set of measure zero where $g$ is ill-behaved.

Staying close to mathematical convention and using everywhere boundedness and ordinary compactness
pays dividends in formalization, allowing us to provide a richer and more robust API leading
to simpler proofs.

Further streamlining is possible by using \lean{fun_prop}.\footnote{\url{https://leanprover-community.github.io/mathlib4_docs/tactics.html\#Mathlib.Meta.FunProp.funPropTacStx}}
This is a \mathlib tactic developed by Tomáš Skřivan
that automatically uses tagged lemmas to show that a function has a certain property
(continuity, differentiability, measurability, etc.).
This tactic is good at writing complicated functions as a composition of simpler functions,
and then proving function properties by proving them first about the simpler functions.

\lean{fun_prop} is not as well-suited for proving integrability
as the other function properties mentioned above,
since the composition and product of integrable functions need not be integrable.
The class \lean{BoundedCompactSupport} is a step towards
making \texttt{fun\_prop} more usable for integrability,
since it is better behaved with respect to multiplication and composition.
After the completion of the Carleson project, some progress was made by Tomáš Skřivan and others to improve \lean{fun_prop} for proving integrability.\href{https://github.com/leanprover-community/mathlib4/pull/39323}{\extlink}

\subsection{Common pitfalls}\label{subsec:pitfalls}

In this project, we encountered several cases of mathematically equivalent representations, and an initially non-obvious choice of the right formalization representation. One such case concerned integrating over a set, i.e.\ with respect to a restricted measure, versus integrating an indicator function.
Both forms are equivalent when the set in question is measurable,\href{https://github.com/leanprover-community/mathlib4/blob/b63493a4746b4651fceb3bcf3a4651cc36a4b8de/Mathlib/MeasureTheory/Integral/Lebesgue/Basic.lean#L496-L503}{\extlink}\href{https://github.com/leanprover-community/mathlib4/blob/b63493a4746b4651fceb3bcf3a4651cc36a4b8de/Mathlib/MeasureTheory/Integral/Bochner/Set.lean#L170-L183}{\extlink}
but the form using a restricted measure is preferred.
The reason is that \lean{s.indicator f}\href{https://github.com/leanprover-community/mathlib4/blob/b63493a4746b4651fceb3bcf3a4651cc36a4b8de/Mathlib/Algebra/Notation/Indicator.lean#L51-L51}{\extlink}
is generally only measurable when the set \lean{s} is measurable,
while if \lean{f} is measurable with respect to the measure \lean{μ}, then it is unconditionally
measurable with respect to \lean{μ.restrict s}.\href{https://github.com/leanprover-community/mathlib4/blob/b63493a4746b4651fceb3bcf3a4651cc36a4b8de/Mathlib/MeasureTheory/Function/StronglyMeasurable/AEStronglyMeasurable.lean#L220-L223}{\extlink}

Similarly, we can represent a periodic function on the real line as \lean{f : ℝ → ℂ} with the property \lean{f.Periodic T}\href{https://github.com/leanprover-community/mathlib4/blob/b63493a4746b4651fceb3bcf3a4651cc36a4b8de/Mathlib/Algebra/Ring/Periodic.lean#L43-L44}{\extlink} and integrate it with respect to the measure \lean{volume.restrict (Ioc 0 T) : Measure ℝ}. Alternatively, we can write it as \lean{f : AddCircle T → ℂ}\href{https://github.com/leanprover-community/mathlib4/blob/b63493a4746b4651fceb3bcf3a4651cc36a4b8de/Mathlib/Topology/Instances/AddCircle/Defs.lean#L187-L188}{\extlink} and integrate it with respect to \lean{volume : Measure (AddCircle T)}.
Depending on the context either of these representations could be more suitable.

Finally, there were several results in the project that required us to consider sets of elements from a given finite type $\texttt{T}$. Since all such sets are finite, from a mathematical point of view there is no difference between the types $\texttt{Set T}$\href{https://github.com/leanprover-community/mathlib4/blob/b63493a4746b4651fceb3bcf3a4651cc36a4b8de/Mathlib/Data/Set/Defs.lean#L44-L51}{\extlink} and $\texttt{Finset T}$.\href{https://github.com/leanprover-community/mathlib4/blob/b63493a4746b4651fceb3bcf3a4651cc36a4b8de/Mathlib/Data/Finset/Defs.lean#L72-L79}{\extlink}
However, these are different types in Lean: $\texttt{Set T}$ is the type of all sets consisting of elements of type $\texttt{T}$, while $\texttt{Finset T}$ is the type of finite sets of elements of $\texttt{T}$. Therefore, we needed to make a decision about which of these types to use in the Carleson project. For the most part, we chose to use $\texttt{Set T}$, although sometimes $\texttt{Finset T}$ is used inside proofs, or when required to apply existing \mathlib results. In hindsight, it would likely have been slightly more convenient to use $\texttt{Finset T}$ throughout, to avoid issues like the following.

If \texttt{T} is endowed with a \texttt{Fintype} instance, then \texttt{Finset.sum} can still be applied to a given \texttt{s : Set T}. There are two possible syntaxes to do this, as the following example shows. Note that the displayed equality is true \emph{propositionally} but not \emph{definitionally}: the left hand side is a sum over \texttt{s.toFinset}, whereas the right hand side is a sum over \texttt{\{p : T | p $\in$ s\}}, and these \texttt{Finset}s are not definitionally equal.
\begin{lstlisting}[mathescape]
example {T : Type} [Fintype T] {s : Set T} [DecidablePred fun p $\mapsto$ p $\in$ s]
	{f : T $\to \mathbb{N}$} : $\sum$ p $\in$ s, f p = $\sum$ p with p $\in$ s, f p
\end{lstlisting}

\section{Future work}\label{sec:future}

After the Carleson project was finished,
we have continued to prove further results in the Carleson repository.
Indeed, \Cref{metric-space-Carleson} is strong enough to prove stronger versions of Carleson's theorem (\Cref{classical-carleson}):
as Hunt \cite{MR238019} showed in 1968,
the assumption in \Cref{classical-carleson} that $f$ is continuous
can be weakened to the statement that $f$ is $L^p$ for $p>1$.
Proving this requires \Cref{real-Carleson} combined with some additional techniques from analysis.

The \emph{Lorentz norm} of a function is
$$\|f\|_{L^{p,q}(X,\mu)}\vcentcolon=
p^{\frac1q}\|t\mu\{|f| > t\}^{\frac1p}\|_{L^q(\mathbb{R}_{>0},\frac{dt}{t})}.$$
It is not hard to see that the $L^{p,p}$-norm corresponds to the $L^p$-norm
and the $L^{p,\infty}$-norm corresponds to the weak $L^p$-norm
(which we also denoted $L^{p,\infty}$ before).
Restricted weak-type estimates such as the one in \Cref{metric-space-Carleson}
specify the boundedness of an operator with respect to the $L^{p,1}$-norm.

There are multiple generalizations of the Marcinkiewicz interpolation theorem that
could be formalized, including multilinear interpolation, abstract interpolation and interpolation between Lorentz spaces.
An interpolation theorem between Lorentz spaces would be especially useful for this project, since we have completed the formalization showing that this result
implies \Cref{classical-carleson} for $L^p$ functions with $p>1$.

Another direction is to formalize further applications of \Cref{linearized-metric-Carleson},
such as the pointwise almost everywhere convergence of the Walsh--Fourier series for functions in $L^2$.

\section{Conclusion}\label{sec:conclusion}

The Carleson project, with a total of 28 contributors,
was one of the largest collaborative formalization efforts of a single mathematical theorem to date,
as well as the first formalization of a state-of-the-art theorem
from the area of harmonic analysis.

Formalizing the metric space Carleson theorem required us
to formalize foundational results from harmonic analysis,
which we are submitting to \mathlib
and which will be useful for the formalization of other results from this field.
Future projects will also benefit from the insights we gained during the project,
including the creation of extended norms
and our reflections about the most convenient ways to work with real number types.

Beyond topic-specific insights,
we also learned lessons from the collaborative nature of this project
which apply to formalization efforts from any area of mathematics.
Our way of assigning and keeping track of ongoing tasks in the Zulip chat was sub-optimal,
and if we were to repeat the project,
we would instead use the existing GitHub issue and project infrastructure (as is done for instance in the PNT$+$ project).
Zulip would still be used as a discussion platform,
and in particular to discuss potential inaccuracies within the blueprint proofs.
To make this more efficient, we would ask one of the topic experts
(in our case, an expert in harmonic analysis)
to frequently monitor this chat and respond to questions,
instead of having F.v.D.\ as an intermediary between the remote formalizers and the harmonic analysis experts.

We would also modify some aspects of the blueprint for the formalization.
We would use \LaTeX\ environments for the main definitions in the project,
to facilitate their discovery,
and we would try to ensure that the results in the blueprint were already proven at the right degree of generality;
this would especially apply to results meant for \mathlib, for which we want a general statement and proof,
rather than a specialization to the project setting.
Finally, we would consider more carefully whether to use non-standard proof methods in the hope of simplifying the formalization,
since this ended up causing most of the mathematical inaccuracies in the first version of the blueprint.

\printbibliography

\end{document}